\documentclass[12pt, a4paper]{amsart}

\usepackage[colorlinks, linkcolor=blue, anchorcolor=blue, citecolor=green]{hyperref}
\usepackage[normalem]{ulem}
\usepackage{graphics,epic}
\usepackage{amsmath,amssymb, amsthm,mathrsfs}
\usepackage[all,2cell]{xy}
\usepackage{shorttoc}
\usepackage{tikz}
\usetikzlibrary{matrix,positioning,decorations.markings,arrows,decorations.pathmorphing,	
	backgrounds,fit,positioning,shapes.symbols,chains,shadings,fadings,calc}
\tikzset{->-/.style={decoration={  markings,  mark=at position #1 with
			{\arrow{>}}},postaction={decorate}}}
\tikzset{-<-/.style={decoration={  markings,  mark=at position #1 with
			{\arrow{<}}},postaction={decorate}}}

\newtheorem{theorem}{Theorem}[section]
\newtheorem*{theorem*}{Theorem}

\newtheorem{lemma}[theorem]{Lemma}
\newtheorem{notation}[theorem]{Notations}
\newtheorem{proposition}[theorem]{Proposition}
\newtheorem{corollary}[theorem]{Corollary}
\newtheorem{conjecture}[theorem]{Conjecture}
\newtheorem*{conjecture*}{Conjecture}

\newtheorem{remark}[theorem]{Remark}
\newtheorem{definition}[theorem]{Definition}

\newtheorem{lem}[theorem]{Lemma}

\newcommand{\opname}[1]{\operatorname{\mathsf{#1}}}

\renewcommand{\mod}{\opname{mod}\nolimits}

\newcommand{\Fac}{\opname{Fac}}

\renewcommand{\P}{\mathbb{P}}

\newcommand{\TT}{\mathbf{T}}
\renewcommand{\SS}{\mathbf{S}}
\newcommand{\MM}{\mathbf{M}}
\newcommand{\ba}{\mathbf{a}}
\newcommand{\PT}{\P(\TT)}
\newcommand{\Hom}{\opname{Hom}}

\newcommand{\Ext}{\opname{Ext}}

\newcommand{\End}{\opname{End}}
\newcommand{\Int}{\opname{Int}}

\newcommand{\mh}{\mathcal{H}}

\newcommand{\nn}{node[black]{$\bullet$}}
\newcommand{\chang}[2]{l_{#1}(#2)}

\newcommand{\Duo}{\operatorname{D}}

\makeatletter
\newcommand*\bigcdot{\mathpalette\bigcdot@{.5}}
\newcommand*\bigcdot@[2]{\mathbin{\vcenter{\hbox{\scalebox{#2}{$\m@th#1\bullet$}}}}}
\makeatother

\begin{document}
	
	\title{On support $\tau$-tilting graphs of gentle algebras}

	\author[Fu]{Changjian Fu}
	\address{Changjian Fu\\Department of Mathematics\\SiChuan University\\610064 Chengdu\\P.R.China}
	\email{changjianfu@scu.edu.cn}
	\author[Geng]{Shengfei Geng}
	\address{Shengfei Geng\\Department of Mathematics\\SiChuan University\\610064 Chengdu\\P.R.China}
	\email{genshengfei@scu.edu.cn}
	\author[Liu]{Pin Liu}
	\address{Pin Liu\\Department of Mathematics\\
		Southwest Jiaotong University\\
		610031 Chengdu \\
		P.R.China}
	\email{pinliu@swjtu.edu.cn}
	\author[Zhou]{Yu Zhou }
	\address{Yu Zhou\\Yau Mathematical Sciences Center\\
		Tsinghua University\\
		100084 Beijing \\
		P.R.China}
	\email{yuzhoumath@gmail.com}
	\subjclass[2010]{}
	\keywords{Gentle algebras; Support $\tau$-tilting modules; Marked surface}
	\maketitle
	
	\begin{abstract}
		Let $A$ be a finite-dimensional gentle algebra over an algebraically closed field. We investigate the combinatorial properties of support $\tau$-tilting graph of $A$. In particular, it is proved that the support $\tau$-tilting graph of $A$ is connected and has the so-called reachable-in-face property. This property was conjectured by Fomin and Zelevinsky for exchange graphs of cluster algebras which was recently confirmed by Cao and Li.
	\end{abstract}
	
	\tableofcontents
	
	\section{Introduction and main results}
	
	\subsection{Background}
	
	Support $\tau$-tilting module is the central notion of $\tau$-tilting theory \cite{AIR}, which generalizes the classical tilting module. This class of modules has been investigated in various contexts \cite{CDT, DF} and has been  found to be deeply connected with other content of representation theory, such as torsion classes, silting objects, $t$-structures (cf. \cite{AIR, BY} for instance).   In contrast to tilting modules, the support $\tau$-tilting module can always be mutated at an arbitrary indecomposable direct summand to obtain a new support $\tau$-tilting module. Therefore, the support $\tau$-tilting modules may have a richer combinatorial structure than tilting modules. In particular, various cluster phenomenon in cluster algebras have been found in representation theory of finite-dimensional algebras via support $\tau$-tilting modules.

	Let $A$ be a finite-dimensional algebra over an algebraically closed field $k$. Denote by $\opname{s\tau-tilt} A$ the set of isomorphism classes of basic support $\tau$-tilting $A$-modules. The  support $\tau$-tilting graph $\mh(\opname{s\tau-tilt} A)$ is the exchange graph of support $\tau$-tilting $A$-modules, which encodes the information of mutations. In particular, the support $\tau$-tilting graph  $\mh(\opname{s\tau-tilt} A)$ is  regular and can be regarded as a counterpart of  the exchange graph of the cluster algebra. One of the basic questions in $\tau$-tilting theory is to determine the number of connected components of $\mh(\opname{s\tau-tilt} A)$. We refer to \cite{AIR, DIJ, BMRRT, FG, QZ,Y}  for recent progress on this question for certain classes of finite-dimensional algebras. All the known results suggest the following folklore conjecture, which strengthens the additive reachability conjecture in cluster algebras (cf. \cite[Remark 5.9]{Qin}).
	\begin{conjecture}\label{c:conj-reach}
		Let $A$ be a connected finite-dimensional $k$-algebra. If there is a path between   $A$ and  $0$ in $\mh(\opname{s\tau-tilt} A)$, then $\mh(\opname{s\tau-tilt} A)$ has precisely one connected component.
		\end{conjecture}
	
	The main purpose of this note is to study combinatorial properties of support $\tau$-tilting graphs for  gentle algebras.  Gentle algebras were introduced by Assem and Skowro\'{n}ski \cite{AS}, as an important class of representation-tame finite-dimensional algebras. The class of gentle algebras is closed under derived equivalence \cite{S, SZ} and has attracted a lot of attention due to its occurrence in various contexts, such as Fukaya categories \cite{HKK}, dimer models \cite{B} and cluster theory \cite{ABCP,BZ}.

	\subsection{Main results}

	Let $k$ be an algebraically closed field and $A$ a finite-dimensional $k$-algebra. Denote by $\mod A$ the category of finitely generated right $A$-modules. Let $\tau$ be the Auslander-Reiten translation of $\mod A$. For any $M\in\mod A$, denote by $|M|$ the number of pairwise non-isomorphic indecomposable summands of $M$; $M$ is called {\it basic} if the number of indecomposable summands of $M$ equals $|M|$; $M$ is called {\it $\tau$-rigid} if $\Hom_A(M,\tau M)=0$; $M$ is called {\it $\tau$-tilting} if it is $\tau$-rigid and $|M|=|A|$. A {\it $\tau$-rigid pair} is a pair $(M,P)$ with $M\in \mod A$ and $P$ a finitely generated projective $A$-module, such that $M$ is $\tau$-rigid and $\Hom_A(P,M)=0$. A $\tau$-rigid pair is called a {\it support $\tau$-tilting pair} if $|M|+|P|=|A|$. In this case, $M$ is a {\it support $\tau$-tilting} $A$-module and $P$ is uniquely determined by $M$ up to isomorphism provided that $P$ is basic. In the following we always identify support $\tau$-tilting modules with support $\tau$-tilting pairs. 
	
	Let $(M,P)$ and $(N,Q)$ be two $\tau$-rigid pairs. We say that $(N,Q)$ is a direct summand of $(M,P)$ if $N$ and $Q$ are direct summands of $M$ and $P$ respectively. A $\tau$-rigid pair $(M,P)$ is {\it indecomposable} if $|M|+|P|=1$. In particular, each basic support $\tau$-tilting pair has $|A|$ indecomposable direct summands. A basic $\tau$-rigid pair $(M,P)$ is {\it almost complete support $\tau$-tilting} if $|M|+|P|=|A|-1$. It has been proved in \cite{AIR} that there exist exactly two non-isomorphic basic support $\tau$-tilting pairs $(M_i,P_i)$ such that $(M,P)$ is a direct summand of $(M_i, P_i)$ for $i=1,2$. In this case, $(M_1, P_1)$ and $(M_2, P_2)$ are {\it mutation} of each other. Clearly, $(M_1, P_1)$ and $(M_2, P_2)$ are only different in one indecomposable direct summand.
	
	\begin{definition}[support $\tau$-tilting graph]
		The {\it support $\tau$-tilting graph} $\mh(\opname{s\tau-tilt} A)$ has vertex set indexed by the isomorphism classes of basic support $\tau$-tilting $A$-modules, and two basic support $\tau$-tilting modules are connected by an edge if and only if they are mutations of each other.
	\end{definition}

	Our first main result is about the number of the connected components of  $\mh(\opname{s\tau-tilt} A)$ of a gentle algebra $A$, which provides new evidences for Conjecture \ref{c:conj-reach}.
	
	\begin{theorem}\label{t:main-thm-1}
		Let $A$ be a finite-dimensional gentle algebra over an algebraically closed field $k$. The support $\tau$-tilting graph $\mh(\opname{s\tau-tilt} A)$ of $A$ has precisely one connected component.
	\end{theorem}
	
	It is known that the graph $\mh(\opname{s\tau-tilt} A)$ is isomorphic to the full subgraph of the silting graph of $A$ consisting of 2-term silting complexes (see \cite{AIR}). Note that the silting quiver of a gentle algebra possibly have infinitely many connected components (see \cite{Du}).

	In \cite{DIJ}, Demonet, Iyama and Jasso introduced a simplicial complex $\Delta(A)$ for a finite-dimensional algebra $A$ via $\tau$-tilting theory. In particular, there is a one-to-one correspondence between the $d$-simplexes of $\Delta(A)$ and the basic $\tau$-rigid pairs of $A$ which have exactly $d+1$ indecomposable summands. The simplicial complex $\Delta(A)$ has pure dimension of $|A|-1$. The support $\tau$-tilting graph $\mh(\opname{s\tau-tilt} A)$ can be identified with the dual graph of $\Delta(A)$. Under this identification, each basic $\tau$-rigid pair $(M,P)$ determines a face $\mathcal{F}_{(M,P)}$ of $\mh(\opname{s\tau-tilt} A)$. Namely, the face $\mathcal{F}_{(M,P)}$ is the full subgraph of $\mh(\opname{s\tau-tilt} A)$ consisting of basic support $\tau$-tilting pairs which admit $(M,P)$ as a direct summand. 
	
	Let $(M,P)$ and $(N,Q)$ be basic support $\tau$-tilting pairs. We say $(N,Q)$ is {\it mutation-reachable} by $(M,P)$ if one can obtain $(N,Q)$ from $(M,P)$ by a finite sequence of mutations. Equivalently, there is a path from $(M,P)$ to $(N,Q)$ in $\mh(\opname{s\tau-tilt} A)$. The following definition is inspired by \cite[Conjecture 4.14(3) ]{FZ} for exchange graphs of cluster algebras, which was recently confirmed in \cite{CL}. We also remark that the definition  can be formulated for an abstract exchange graph in \cite{BY}.
	
	\begin{definition}[reachable-in-face]
		The support $\tau$-tilting graph $\mh(\opname{s\tau-tilt} A)$ has the {\it reachable-in-face} property if for any mutation-reachable basic support $\tau$-tilting pairs $(M,P)$ and $(N,Q)$ such that $(M,P)$ and $(N,Q)$ have a common direct summand $(L,R)$, there is a path from $(M,P)$ to $(N,Q)$ lying in the face $\mathcal{F}_{(L,R)}$ determined by $(L,R)$.  In this case, we also say that $A$ has the reachable-in-face property.
	\end{definition}
	The definition of reachable-in-face property is related to the non-leaving-face property, which was introduced in \cite{CP} for polytopes. Here we generalize it to the support $\tau$-tilting graph of an arbitrary finite-dimensional $k$-algebra $A$. We say that  the support $\tau$-tilting graph $\mh(\opname{s\tau-tilt} A)$ has the {\it non-leaving-face} property provided that for any mutation-reachable basic support $\tau$-tilting pairs $(M,P)$ and $(N,Q)$ with maximal common direct summand $(L,R)$, every path with minimal length linked $(M,P)$ and $(N,Q)$ lies in the face $\mathcal{F}_{(L,R)}$ determined by $(L,R)$. Clearly, the non-leaving-face property implies the reachable-in-face property.  It is worth mentioning that  for a $2$-Calabi-Yau tilted gentle algebra $A$ arising from a marked surface without punctures, Br\"{u}stle and Zhang \cite{BZ2} have proved the non-leaving-face property for the support $\tau$-tilting graph of $A$.

	For an arbitrary finite-dimensional $k$-algebra, we do not know whether its support $\tau$-tilting graph has the reachable-in-face property. We do know that the reachable-in-face property holds for the following class of algebras:
	\begin{enumerate}
		\item $\tau$-tilting finite algebras \cite{DIJ}, i.e., algebras with finitely many pairwise non-isomorphic indecomposable $\tau$-rigid modules;
		\item Cluster-tilted algebras arising from hereditary abelian categories \cite{FG};
		\item More generally, $2$-Calabi-Yau tilted algebras \cite{C}.
	\end{enumerate}
	
	Our second result shows that the reachable-in-face property holds for gentle algebras.
	\begin{theorem}\label{t:main-thm-2}
		Let $A$ be a finite-dimensional gentle algebra over an algebraically closed field $k$.  The support $\tau$-tilting graph $\mh(\opname{s\tau-tilt} A)$ has the reachable-in-face property.
	\end{theorem}

	Our proofs of Theorem \ref{t:main-thm-1} and Theorem \ref{t:main-thm-2} are inspired by the reduction approach in \cite{FG}, where the Iyama-Yoshino's reduction was applied to study the connectedness of  cluster-tilting graph of a hereditary category. In present paper,  the Iyama-Yoshino's reduction has been replaced by $\tau$-reduction in the sense of Jasso \cite{J}. We introduce the notion of $\tau$-reachable and totally $\tau$-reachable for an arbitrary finite-dimensional $k$-algebra (cf. Definition \ref{d:tau-reachable} and Definition \ref{d:totally-tau-reachable}) and prove Theorem \ref{t:main-thm-1} and Theorem \ref{t:main-thm-2} by showing that every finite-dimensional gentle $k$-algebra has the totally $\tau$-reachable property.  In particular, we have the following general result.

	\begin{theorem}[Theorem \ref{t:totally-equivalent-to-connected}]\label{t:main-thm-3}
		Let $A$ be a finite-dimensional $k$-algebra. Then $A$ is totally $\tau$-reachable if and only if $A$ has the reachable-in-face property and the support $\tau$-tilting graph $\mh(\opname{s\tau-tilt} A)$ is connected.
	\end{theorem}

	The paper is organized as follows. In Section \ref{s:gg}, we recall basic properties of gentle algebras and their geometric models. In Section \ref{s:reachability-property}, we introduce the definition of $\tau$-reachability for indecomposable $\tau$-rigid pairs and establish the $\tau$-reachable property for gentle algebras. In Section \ref{s:proofs-main-results}, we first give the proof of Theorem \ref{t:main-thm-3} and then deduce Theorem \ref{t:main-thm-1} and Theorem \ref{t:main-thm-2}.
	The reduction theory of $\tau$-rigid pairs is presented in Appendix \ref{s:reduction}.
	
	\subsection*{Acknowledgments}
	The authors thank Professor Wen Chang and Professor Thomas Br\"{u}stle  for their interest and helpful comments. They are grateful to H\aa vard Utne Terland for comments on the converse of Proposition \ref{p:connect-reachable}. This work is partially supported by the National Natural Science Foundation of China (Grant No. 11801297, 11971326, 12071315).

	\section{Gentle algebras and their geometric models}\label{s:gg}
	Let $Q$ be a quiver. Denote by $Q_0$  the set of vertices and $Q_1$ the set of arrows.  For an arrow $\alpha\in Q_1$, we denote by $s(\alpha)$ and $t(\alpha)\in Q_0$ the source and target of $\alpha$ respectively.
	
	\subsection{Gentle algebras}
	A finite-dimensional $k$-algebra $A$ is {\it gentle} if it admits a presentation $A=kQ/I$ satisfying the following conditions:
	\begin{itemize} 
		\item[(G1)] Each vertex of $Q$ is the source of at most two arrows and the target of at most two arrows;
		\item[(G2)] For each arrow $\alpha$, there is at most one arrow $\beta$ (resp. $\gamma$) such that $\beta\alpha\in I$ (resp. $\gamma\alpha\not\in I$);
		\item[(G3)] For each arrow $\alpha$, there is at most one arrow $\beta$ (resp. $\gamma$) such that $\alpha\beta\in I$ (resp. $\alpha\gamma\not\in I$);
		\item[(G4)] I is generated by paths of length 2.
	\end{itemize}
	
	
	The following is a direct consequence of definition of gentle algebras.
	\begin{lem}\label{l:factor-gentle-algebra}
		Let $A$ be a finite-dimensional gentle algebra over $k$ and $e$ an idempotent of $A$. The factor algebra $A/\langle e\rangle$ is also gentle.
	\end{lem}
	
	The following result plays a fundamental role in our reduction approach.
	
	\begin{lemma}[{\cite[Theorem 1.1]{S}}]\label{l:endomorphism algebra is gentle}
		Let $A$ be a finite-dimensional gentle algebra over $k$ and ${M}$ a rigid $A$-module (\text{i.e.}, $\Ext^1_A(M,M)=0$), the endomorphism algebra $\End_A (M)$ is a gentle algebra.
	\end{lemma}

	\subsection{Tilings and permissible curves}\label{s:gentle-algebra}
	
	In this subsection, we mainly follow \cite{BS} and \cite{HZZ}.
	
	A {\it marked surface} is a pair $(\SS, \MM)$, where $\SS$ is a connected oriented Riemann surface with non-empty boundary $\partial \SS$, and $\MM\subset\partial\SS$ is a finite set of marked points on the boundary. A connected component of $\partial \SS$ is called a {\it boundary component} of $\SS$. A boundary component $B$ of $\SS$ is {\it unmarked} if $\MM\cap B=\emptyset$. We allow unmarked boundary components of $\SS$. A {\it boundary segment} is the closure of a component of $\partial\SS\backslash\MM$.
	
	An {\it arc} on $(\SS, \MM)$ is a continuous map $\gamma: [0,1]\to \SS$ such that
	\begin{itemize}
		\item[$\circ$] $\gamma(0),\gamma(1)\in \MM$ and $\gamma(t)\in \SS\backslash \MM$ for $0<t<1$;
		\item[$\circ$] $\gamma$ is neither null-homotopic nor homotopic to a boundary segment.
	\end{itemize}
	We always consider arcs on $\SS$ up to homotopy relative to their endpoints and up to inverse, where the {\it inverse} of an arc $\gamma$ on $\SS$ is defined as $\gamma^{-1}(t)=\gamma(1-t)$ for $t\in [0,1]$. Denote by $\mathbf{C}(\SS)$ the set of arcs on $\SS$. An arc $\gamma$ is called a {\it loop} if $\gamma(0)=\gamma(1)$.

	For any arcs $\gamma_1, \gamma_2\in\mathbf{C}(\SS)$ which are in a minimal position, the {\it intersection number} between them is defined to be
	\[\Int(\gamma_1,\gamma_2):=|\{(t_1,t_2)~|~0<t_1,t_2<1, \gamma_1(t_1)=\gamma_2(t_2)\}|.
	\]
	An arc $\gamma\in\mathbf{C}(\SS)$ is said to be without self-intersections if $\Int(\gamma,\gamma)=0$.
	
\begin{figure}[htpb]
		\begin{tikzpicture}

			\draw[thick,fill=black!20] (2,0.17)arc (250:290:3);
			
			\fill(3,0) circle(1.5pt);

			\draw[thick,color=blue] (3,0) .. controls (-1,-3) and (8,-3) .. (3,0);
			
			\draw[thick,fill=black!20] (3,-1.3) circle(8pt);
			\node at (3,-3){Type I};

			\draw[thick,fill=black!20] (9,0.17)arc (250:290:3);
			
			\fill(10,0) circle(1.5pt);
			
			\draw[thick,fill=black!20] (9,-2.17)arc (110:70:3);
			\fill(10,-2) circle(1.5pt);
			\draw[thick,fill=black!20] (10,-1) circle(8pt);
			\draw[thick,color=blue] (10,0) .. controls (9,-1) and (9.5,-1.5) .. (10,-2);
			\draw[thick,color=blue] (10,0) .. controls (11.5,-1) and (11,-1.5) .. (10,-2);
			\node at (10,-3){Type II};

		\end{tikzpicture}
		\caption{Tiles of type I and II}\label{f:tile-type}
	\end{figure}
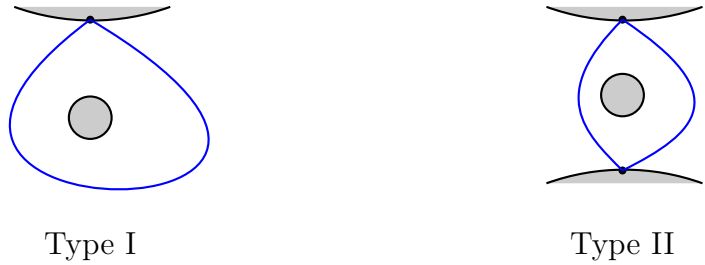

	A {\it partial triangulation} $\TT$ of $(\SS, \MM)$ is a collection of arcs in $\mathbf{C}(\SS)$ without self-intersections and such that $\Int(\gamma_1,\gamma_2)=0$ for any $\gamma_1,\gamma_2\in\TT$. A triple $(\SS, \MM, \TT)$ is called a {\it tiling} provided that $(\SS,\MM)$ is a marked surface and $\TT$ is a partial triangulation of $(\SS, \MM)$ such that $\SS$ is divided by $\TT$ into a collection of regions (also called {\it tiles}) of the following types:
	\begin{enumerate}
		\item[(I)] monogons with exactly one unmarked boundary component in their interior (see the left picture in Figure \ref{f:tile-type});
		\item[(II)] digons with exactly one unmarked boundary component in their interior (see the right picture in Figure \ref{f:tile-type});
		\item[(III)] $m$-gons, with $m\ge 3$ and whose edges are arcs of $\TT$ and at most one boundary segment, and whose interior contains no unmarked boundary component of $\SS$;
		\item[(IV)] $3$-gons bounded by two boundary segments and one arc of $\TT$, and whose interior contains no unmarked boundary component of $\SS$.
	\end{enumerate}

	An {\it arc segment} in a tile $\Duo$ is a curve $\eta:[0,1]\to\Duo$ such that $\eta(0)$ and $\eta(1)$ are in the edges of $\Duo$, and $\eta(t),0<t<1$ are in the interior of $\Duo$. An arc segment in $\Duo$ is called {\it permissible} with respect to $\TT$ if it satisfies one of the following conditions (cf. \cite[Definition 2.1]{HZZ} and compare \cite[Definition~3.1]{BS}).
	\begin{enumerate}
		\item[(P1)] One endpoint $P$ of $\eta$ is in $\MM$ and the other $Q$ is in the interior of a non-boundary edge, say $\gamma$ of $\Duo$, such that $\eta$ is not isotopic to a segment of an edge of $\Duo$ relative to their endpoints, and after moving $P$ along the edges of $\Duo$ in anticlockwise order to the next marked point, say $P'$, the new arc segment obtained from $\eta$ is isotopic to a segment of $\gamma$ relative their endpoints. See Figure~\ref{f:P1} for all the possible cases of permissible arc segments satisfying (P1).
		
		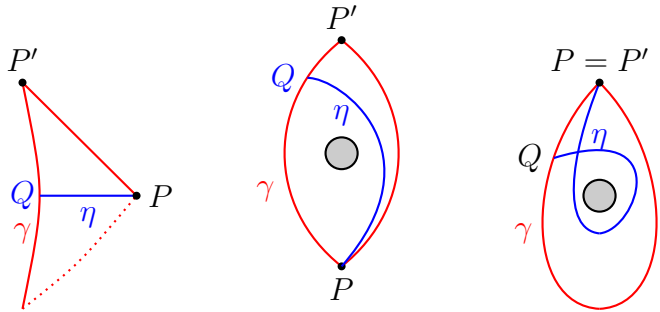
\begin{figure}
			\begin{tikzpicture}
				\draw[thick,color=red] (1,0) .. controls (1.3,1.5) and (1.3,1.5) .. (1,3);
				\draw[thick,color=red] (1,3)node[black,above]{$P'$} .. controls (2,2) and (2, 2) .. (2.5,1.5);
				\draw[thick,color=red,dotted] (1,0) .. controls (1.75, 0.6) and (1.75,0.6) .. (2.5,1.5);
				\draw[thick,color=blue] (2.5,1.5)node[black,right]{$P$} tonode[below]{$\eta$} (1.23,1.5);
				\fill(2.5,1.5) circle(1.5pt); \fill(1,3) circle(1.5pt);
				\draw[red] (1,1)node{$\gamma$} (1,1.5)node[blue]{$Q$};	
			\end{tikzpicture}\qquad
			\begin{tikzpicture}
				\draw[thick,color=red] (7,0) .. controls (8,0.8) and (8,2.2) .. (7,3)node[black,above]{$P'$};	
				\draw[thick,color=red] (7,0) .. controls (6,0.8) and (6,2.2) .. (7,3);	
				\fill(7,0) circle(1.5pt);		\fill(7,3) circle(1.5pt);	
				\draw[thick,fill=black!20] (7,1.5) circle(6pt);
				\draw[thick,color=blue] (7,0)node[black,below]{$P$} .. controls (8.3,1.5) and (7,2.5) .. (6.55,2.5)node[blue,left]{$Q$};
				\draw[blue] (7,2)node{$\eta$};
				\draw[red] (6,1)node{$\gamma$};
			\end{tikzpicture}\qquad
			\begin{tikzpicture}
				\draw[thick,color=red] (13,0) .. controls (14,0.05) and (14,2) .. (13,3)node[black,above]{$P=P'$};
				\draw[thick,color=red] (13,0) .. controls (12,0.05) and (12,2) .. (13,3);
				\draw[thick,color=blue] (13,1) .. controls (12.5,1.05) and (12.6,2) .. (13,3);	
				\draw[thick,color=blue] (13,1) .. controls (13.5,1.05) and (14,2.5) .. (12.4,2)node[black,left]{$Q$};	
				\draw[thick,fill=black!20] (13,1.5) circle(6pt);
				\fill(13,3) circle(1.5pt);
				\draw[blue] (13,2.25)node{$\eta$};
				\draw[red] (12,1)node{$\gamma$};
			\end{tikzpicture}
			\caption{Condition $(P1)$}\label{f:P1}
		\end{figure}
		
		\item[(P2)] The endpoints of $\eta$ are in the interiors of non-boundary edges $x,y$ (which are possibly not distinct) of $\Duo$ such that $\eta$ has no self-intersections, $x,y$ have a common endpoint $p_\eta\in \MM$ and $\eta$ cuts out an angle from $\Duo$ as shown in Figure \ref{f:PAS}. We denote by $\triangle(\eta)$ the local triangle cut out by $\eta$.
		\begin{figure}[htpb]
			\begin{tikzpicture}
				\draw[thick,fill=black!20] (9,0.17)arc (250:290:3);
				\node at (10,0.2){$p_{\eta}$};
				\fill(10,0) circle(1.5pt);
				\draw[thick,color=red] (10,0) .. controls (8.6,-0.9) and (8.4,-1.1) .. (8,-2.2);
				\draw[thick,color=red,dashed] (8,-2.2) -- (7.8,-3);
				\draw[thick,color=red] (10,0) .. controls (10.4,-0.9) and (10.6,-1.1) .. (11,-2.2);
				\draw[thick,color=red,dashed] (11,-2.2) arc (20:-40:1);
				\draw[thick,color=blue] (8,-1) .. controls (9,-1.5) and (10,-1) .. (12,-1);
				\draw[thick,color=blue,dashed] (12,-1) arc (90:50:1.5);
				\draw[thick,color=blue,dashed] (8,-1) -- (7,-0.5);
				\node at (9.5, -1.5){$\eta$};
				\node at (7.5,-3){{\color{red}$x$}};
				\node at (10.5, -3){{\color{red}$y$}};
				\draw (9.7,-.8)node{$\triangle(\eta)$};
			\end{tikzpicture}
			\caption{Condition (P2)}\label{f:PAS}
		\end{figure}
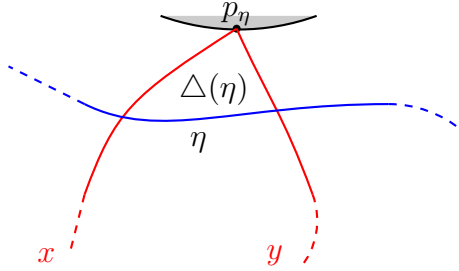
	\end{enumerate}

	For any arc $\gamma\in\mathbf{C}(\SS)$, we always assume that $\gamma$ is in a minimal position with $\TT$.
	\begin{definition}[{\cite[Definition 3.1]{BS} and \cite[Definition 2.6]{HZZ}}]
		An arc $\gamma$ on $\SS$ is called {\it permissible} (with respect to $\TT$) if each arc segment of $\gamma$ divided by $\TT$ is permissible.
	\end{definition}
	
	We denote by $\PT\subset\mathbf{C}(\SS)$ the set of permissible arcs on $\SS$.  For any $\gamma\in\PT$ and any $\ba\in\TT$, we define $\chang{\ba}{\gamma}=\Int(\gamma,\ba)$. The {\it length} of $\gamma$ (with respect to  $\TT$) is defined to be $$\chang{}{\gamma}=\sum_{\ba\in\TT}\chang{\ba}{\gamma}.$$
	
	\subsection{Geometric interpretation of $\tau$-tilting theory}

	
	
	For any finite-dimensional $k$-algebra $A$, denote by $\opname{\tau-}\opname{rigid} A$ the set of isomorphism classes of basic $\tau$-rigid modules in $\mod A$, and by $\opname{ind\tau-}\opname{rigid} A$ the subset of $\opname{\tau-}\opname{rigid} A$ consisting of indecomposable ones.
	
	The following result is useful to study $\tau$-rigid modules of gentle algebras.
	
	\begin{lemma}\label{l:intersection-hom}
		Let $A=k Q/I$ be a finite-dimensional gentle algebra. Then there exists a tiling $(\SS,\MM,\TT)$ such that 
		\begin{itemize}
			\item we have a complete set of primitive orthogonal idempotents $\{e_{\ba}\mid \ba\in\TT\}$ of $A$ indexed by $\TT$, and
			\item there is a bijection
			$$M\colon\{\gamma\in\PT\mid\Int(\gamma,\gamma)=0\}\to\opname{ind\tau-}\opname{rigid} A,$$
			satisfying $\chang{\ba}{\gamma}=\dim\Hom_{A}(e_{\ba}A,M(\gamma))$ for any $\ba\in\TT$ and $\gamma\in\PT$.
		\end{itemize} In particular, $\dim M(\gamma)=\chang{}{\gamma}$ for any $\gamma\in\PT$. Moreover, the bijection $M$ induces a bijection
		$$\{R\subset\PT\mid\Int(\gamma_1,\gamma_2)=0,\ \forall\gamma_1,\gamma_2\in R\}\to\opname{\tau-}\opname{rigid} A,$$
		mapping $R$ to $\oplus_{\gamma\in R}M(\gamma)$.
	\end{lemma}
	
	\begin{proof}
		By \cite[Theorems 2.10 and 3.8]{BS}, there is a tiling $(\SS, \MM,\TT)$ with $|\TT|=|Q_0|$ and such that to any $\gamma\in \PT$, there is associated indecomposable module $M(\gamma)\in\mod A$ satisfying $\Int(\gamma,\ba)=\dim\Hom_{A}(e_{\ba}A,M(\gamma))$ for any $\ba\in\TT$. The two bijections in the lemma then follows from \cite[Propositions 5.3 and 5.6]{HZZ}.
	\end{proof}

	\begin{remark}
		In \cite[Theorem 2.10]{BS}, any finite-dimensional gentle algebra is realized as a tiling algebra; in \cite[Theorem 3.8]{BS}, $M$ is a bijection from $\PT$ to the set of non-zero string modules in $\mod A$; the bijections in \cite[Propositions 5.3 and 5.6]{HZZ} also exist for skew-gentle algebras.
	\end{remark}

	We fix some notations for a permissible arc in a tiling.
	
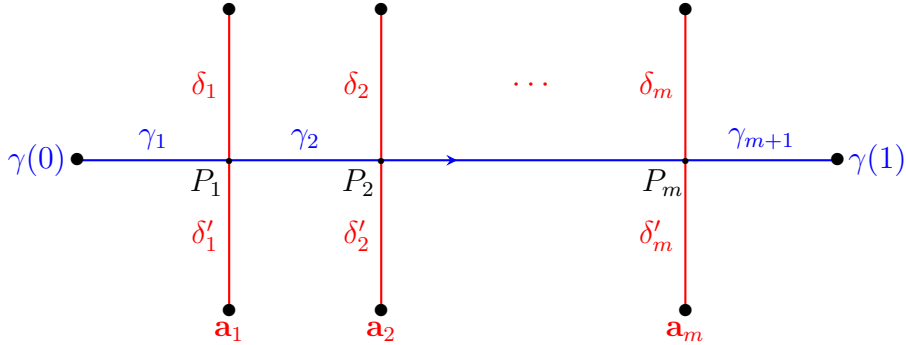
\begin{figure}[htpb]
		\begin{tikzpicture}
			\draw[blue,thick,->-=.5,>=stealth] (-5,0)\nn to (5,0)\nn;
			\draw[blue] (-5,0)node[left]{$\gamma(0)$} (5,0)node[right]{$\gamma(1)$};
			\draw[red,thick] (-3,-2)\nn to (-3,2)\nn;
			\draw[red,thick] (-1,-2)\nn to (-1,2)\nn;
			\draw[red] (1,1)node{$\cdots$};
			\draw[red,thick] (3,-2)\nn to (3,2)\nn;
			\draw[red] (-3,-2)node[below]{$\ba_1$} (-1,-2)node[below]{$\ba_2$} (3,-2)node[below]{$\ba_m$};
			\draw[blue] (-4,0)node[above]{$\gamma_1$} (-2,0)node[above]{$\gamma_2$} (4,0)node[above]{$\gamma_{m+1}$};
			\draw (-3.3,-.3)node{$P_1$} (-1.3,-.3)node{$P_2$} (2.7,-.3)node{$P_m$};
			\draw[red] (-3,1)node[left]{$\delta_1$} (-1,1)node[left]{$\delta_2$} (3,1)node[left]{$\delta_m$};
			\draw[red] (-3,-1)node[left]{$\delta_1'$} (-1,-1)node[left]{$\delta_2'$} (3,-1)node[left]{$\delta_m'$};
			\draw (-3,-.05)node{$\bigcdot$} (-1,-.05)node{$\bigcdot$} (3,-.05)node{$\bigcdot$};
		\end{tikzpicture}
		\caption{Notations for a permissible arc $\gamma$}\label{f:notation}
	\end{figure}

	\begin{notation}\label{n:seg}
		Let $(\SS,\MM,\TT)$ be a tiling and $\gamma\in\PT$. Let $m=\chang{}{\gamma}$. Fix an orientation of $\gamma$. We denote by $P_1,\ldots,P_m$ the intersections of $\gamma$ and $\TT$ in order, which are in $\ba_1,\ldots,\ba_m\in\TT$ respectively, and which divide $\gamma$ into arc segments $\gamma_1, \ldots, \gamma_{m+1}$ in order. Note that we may have $\ba_i=\ba_j$ for different $i$ and $j$. 
		
		For any $1\leq j\leq m$, we denote by $\delta_j$ (resp. $\delta_j'$) the left half (resp. the right half) segment of $\ba_j$ divided by $P_j$, where left/right is w.r.t. the orientation of $\gamma$. See Figure~\ref{f:notation}. We fix the orientations of each $\delta_j,\delta_j'$ such that $P_j=\delta_j(1)=\delta_j'(1)$.
	\end{notation}
	
	Using the notations in Notations~\ref{n:seg}, we call $P_j$ is {\it left most} if $\delta_j$ does not contain any other $P_i$.
	
	We simply denoted by $\gamma_1\circ \gamma_2$ the {\it concatenation} of $\gamma_1:[0,1]\to\SS$ and $\gamma_2:[0,1]\to\SS$ with $\gamma_1(1)=\gamma_2(0)$. The following lemma is useful to construct new permissible arcs from a given one, which will be used in the next section.
	
	\begin{lemma}\label{l:concatenation}
		Let $\gamma\in\PT$ without self-intersections. Using the notations in Notations~\ref{n:seg}, for any $1<j<m+1$, either $\delta_{j-1}$, $\delta_{j}$ and $\gamma_j$, or $\delta_{j-1}'$, $\delta_{j}'$ and $\gamma_j$, are the three edges of the triangle $\triangle(\gamma_j)$.
		\begin{itemize}
			\item[(a)] In the former case, if $P_{j-1}$ is left most, then so is $P_{j}$.
			\item[(b)] In the later case, the concatenation $\delta_{j-1}\circ\gamma_j$ is isotopic to a permissible arc segment satisfying condition $(P1)$.
		\end{itemize}
	\end{lemma}
	\begin{proof}
		Since $\gamma_j$ is a permissible arc segment satisfying condition (P2), it cuts out an angle formed by either $\delta_{j-1}$ and $\delta_{j}$, or $\delta_{j-1}'$ and $\delta_{j}'$. So the first assertion holds, see Figure~\ref{fig:proof-concat}. 
		
		\begin{figure}[htpb]
			\begin{tikzpicture}
				\draw[red,thick] (1.5,-2)\nn to (0,2)\nn to (-1.5,-2)\nn;
				\draw[red] (-1.5,-2)node[below]{$\ba_{j-1}$} (1.5,-2)node[below]{$\ba_{j}$};
				\draw[blue,thick,->-=.5,>=stealth] (-2,0)to(2,0);
				\draw[blue,thick,dashed] (-3,0)node[left]{$\gamma$}to(-2,0) (2,0)to(3,0);
				\draw (-1.3,-.3)node{$P_{j-1}$} (1.1,-.3)node{$P_{j}$};
				\draw (0,.5)node{$\triangle(\gamma_j)$};
				\draw[blue] (0,0)node[below]{$\gamma_j$};
				\draw[red] (-.9,1)node{$\delta_{j-1}$} (.8,1)node{$\delta_{j}$};
				\draw[red] (-1.5,-1)node{$\delta_{j-1}'$} (1.4,-1)node{$\delta_{j}'$};
				\draw (-.75,-.05)node{$\bigcdot$} (.75,-.05)node{$\bigcdot$};
			\end{tikzpicture}\qquad
			\begin{tikzpicture}
				\draw[red,thick] (1.5,2)\nn to (0,-2)\nn to (-1.5,2)\nn;
				\draw[red] (-1.5,2)node[above]{$\ba_{j-1}$} (1.5,2)node[above]{$\ba_{j}$};
				\draw[blue,thick,->-=.5,>=stealth] (-2,0)to(2,0);
				\draw[blue,thick,dashed] (-3,0)node[left]{$\gamma$}to(-2,0) (2,0)to(3,0);
				\draw (-1.3,.3)node{$P_{j-1}$} (1.1,.3)node{$P_{j}$};
				\draw (0,-.5)node{$\triangle(\gamma_j)$};
				\draw[blue] (0,0)node[above]{$\gamma_j$};
				\draw[red] (-.9,-1)node{$\delta_{j-1}'$} (.8,-1)node{$\delta_{j}'$};
				\draw[red] (-1.5,1)node{$\delta_{j-1}$} (1.4,1)node{$\delta_{j}$};
				\draw (-.75,-.05)node{$\bigcdot$} (.75,-.05)node{$\bigcdot$};
			\end{tikzpicture}
			\caption{Cases in Lemma \ref{l:concatenation}}\label{fig:proof-concat}
		\end{figure}
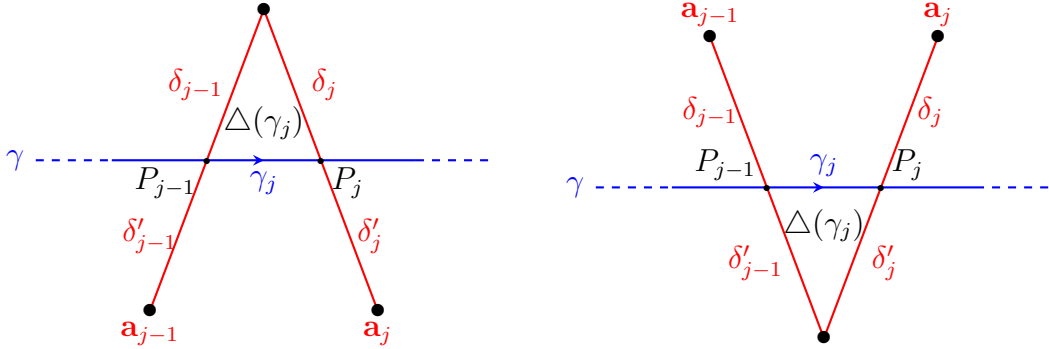
		
		In the former case, since $\delta_{j-1}$, $\delta_{j}$ and $\gamma_j$ are the edges of $\triangle(\gamma_j)$ in the clockwise order, $\delta_{j-1}\circ\gamma_j$ is isotopic to $\delta_j$. Suppose that $P_{j-1}$ is left most. If the interior of $\delta_j$ contains a $P_i$, then there is an arc segment $\eta$ of $\gamma$ in the interior of $\triangle(\gamma_j)$ which has $P_i$ as an endpoint. Since $\gamma$ does not have self-intersection, $\eta$ does not cross $\gamma_j$. Then the other endpoint of $\eta$ has to be in the interior of $\delta_{j-1}$, a contradiction. So $P_j$ is also left most.
		
		In the latter case, since the concatenation $\ba_{j-1}\circ(\delta_{j-1}\circ\gamma_j)$ is isotopic to $\delta'_j$, where the orientation of $\ba_{j-1}$ is taken to be the direction from $\delta_{j-1}'$ to $\delta_{j-1}$. By definition, $\delta_{j-1}\circ\gamma_j$ is a permissible arc segment satisfying condition (P1). 
	\end{proof}
	
	\section{$\tau$-reachable property for gentle algebras}\label{s:reachability-property}
	\subsection{Definition of $\tau$-reachability}
	Let $A$ be a finite-dimensional $k$-algebra. Denote by $\mod A$ the category of finitely generated right $A$-modules.
	
	\begin{definition}
		Let $M$ and $N$ be $\tau$-rigid $A$-modules. We say that $M$ is {\it $\tau$-reachable} from $N$, denoted by $\xymatrix{M\ar@{~}[r]^{\tau} &N,}$ if there exists a sequence of $\tau$-rigid $A$-modules $M_1,\ldots, M_s$ such that $M\oplus M_1, M_1\oplus M_2,\ldots, M_{s-1}\oplus M_s, M_s\oplus N$ are $\tau$-rigid  in $\mod A$, where $s\geq 0$.  
	\end{definition}
	By definition, $M$ is $\tau$-reachable from $N$ if and only if $N$ is $\tau$-reachable from $M$.  All direct summands of a $\tau$-rigid module are $\tau$-reachable from each other. We extend the definition to $\tau$-rigid pairs.
	
	\begin{definition}\label{d:tau-reachable}
		Let $(M,P),(N,Q)$ be two $\tau$-rigid pairs in $\mod A$. We say that $(M,P)$ is {\it $\tau$-reachable} from $(N,Q)$, denoted by $\xymatrix{(M,P)\ar@{~}[r]^{\tau}&(N,Q),}$ if there exists a sequence of $\tau$-rigid pairs $(M_1,P_1),\ldots, (M_s,P_s)$ such that $(M\oplus M_1, P\oplus P_1), (M_1\oplus M_2, P_1\oplus P_2), \ldots, (M_{s-1}\oplus M_s,P_{s-1}\oplus P_s), (M_s\oplus N, P_s\oplus Q)$ are $\tau$-rigid  pairs in $\mod A$.
	\end{definition}
	
	Let $(M,P)$ be a $\tau$-rigid pair. By definition, $(M,0)$ is $\tau$-reachable from $(0,P)$. Let $M, N$ be $\tau$-rigid modules, then $\xymatrix{M\ar@{~}[r]^{\tau} &N}$ implies  $\xymatrix{(M,0)\ar@{~}[r]^{\tau} &(N,0).}$  It is also clear that $\tau$-reachability is reflexive, symmetric and transitive.
	\begin{definition}
		The algebra $A$ has the {\it $\tau$-reachable property} if any two indecomposable $\tau$-rigid pairs of $\mod A$ are $\tau$-reachable from each other.
	\end{definition}
	
	The above definitions are inspired by \cite[Definition 4.1]{FG}. The $\tau$-reachable property has a close relation with the connectedness of support $\tau$-tilting graph. In particular, by definition, we have
	
	\begin{proposition}\label{p:connect-reachable}
		Let $A$ be a finite-dimensional $k$-algebra such that $|A|>1$. If the support $\tau$-tilting graph of $A$ is connected, then $A$ has the $\tau$-reachable property.
	\end{proposition}

	In particular, every $\tau$-tilting finite algebra $A$ with $|A|>1$ has the $\tau$-reachable property. 
		\begin{remark}\label{rem:converse-reachable}
			For connected algebra, we don’t know whether the converse of Proposition 3.4 is true, but for disconnected algebra, the converse is false. In fact, suppose that $A=A_1 \times A_2$, where $A_1$ and $A_2$ are finite dimensional $k$-algebras. For each indecomposable $\tau$-rigid $A_1$-module $M_1$, the pair $(M_1, 0)$ is $\tau$-reachable from $(0, A_2) $ since $(M_1, 0),(M_1, A_2), (0, A_2)$ is exactly the desired sequence. Similarly, for each indecomposable $\tau$-rigid $A_2$-module $M_2$, the pair $(M_2, 0)$ is $\tau$-reachable from $(0, A_1)$. Thus the disconnected finite-dimensional $k$-algebra $A$ has the $\tau$-reachable property. On the other hand, as stated in \cite[Lemma 2.6]{T}, the support $\tau$-tilting graph of the disconnected algebra $A$ has $n_1n_2$ connected components, where $n_1$ (resp. $n_2$) is the number of connected components of $\mh(\opname{s\tau-tilt} A_1)$ (resp. $\mh(\opname{s\tau-tilt} A_2)$).
		\end{remark}

	In the next subsection, we will establish the $\tau$-reachable property for gentle algebras, which plays a key role in the proof of the connectedness of support $\tau$-tilting graphs of gentle algebras. We remark that there exist algebras which do not satisfy the $\tau$-reachable property. 
	
	\begin{remark}\label{rmk:rank1}
		If $|A|=1$, then the $\tau$-rigid pair $(0,0)$ is almost complete. By \cite[Theorem~0.4]{AIR}, there are exactly two non-isomorphic support $\tau$-tilting modules $A$ and $0$, which are related by a mutation. Hence the support $\tau$-tilting graph is connected. However, by definition, $A$ does not have  the $\tau$-reachable property.
	\end{remark}

	\subsection{$\tau$-reachability of gentle algebras}
	Throughout this subsection, let $A$ be a finite-dimensional gentle algebra over $k$. Recall that a module $M$ is {\it sincere } if  $\Hom(P, M)\neq 0$ for any indecomposable projective $A$-module $P$.
	
	\begin{lemma}\label{l:sincere-reachable}
		Assume that  $|A|>1$. Let $M$ be an indecomposable sincere $\tau$-rigid $A$-module. Then there is an indecomposable $\tau$-rigid $A$-module $N$ with $\dim N<\dim M$  such that $(M,0)$ is  $\tau$-reachable from $(N,0)$.
	\end{lemma}
	\begin{proof}
		Let $(\SS, \MM, \TT)$ be the tiling for $A$ given in Lemma~\ref{l:intersection-hom}. Then there is a permissible arc $\gamma$ without self-intersections and such that $M\cong M(\gamma)$. Let $m=\chang{}{\gamma}=\dim M$. Since $M$ is sincere, we have $\Int(\gamma, \ba)\neq 0$ for any $\ba\in \TT$. By the assumption $|A|>1$, we have $|\TT|>1$ and $m>1$. We are going to construct a permissible arc $\gamma'$ without self-intersections such that $\Int(\gamma, \gamma')=0$ and $\chang{}{\gamma'}<\chang{}{\gamma}$, or two permissible arcs $\gamma',\gamma''$ without self-intersections such that $\Int(\gamma,\gamma')=0$, $\chang{}{\gamma'}=\chang{}{\gamma}$ and $\Int(\gamma'',\gamma')=0$, $\chang{}{\gamma''}<\chang{}{\gamma'}$.
		
		Fixing an orientation of $\gamma$ and using the notations in Notation~\ref{n:seg}, we divide the proof into several cases.

		\noindent{\bf Case 1: there is an arc $\ba\in\TT$ which contains at least one of $P_1$ and $P_m$ and at least one of $P_j,1<j<m$.} 
		
		Reversing the orientation of $\gamma$ if necessary, we may assume $\ba$ contains $P_1$ and $P_j$ for some $1<j<m$, such that there is no $P_i$ in $\ba$ between $P_1$ and $P_j$. See Figure~\ref{f:subcase1.1}. Denote by $\delta$ the segment of $\ba$ from $P_1$ to $P_j$. Let $g$ be the segment of $\gamma$ from $P_j$ to one endpoint of $\gamma$ such that the arc segment of $g$ connecting $P_j$ is on the different side of $\ba$ from $\gamma_1$. Take $\gamma'=\gamma_1\circ\delta\circ g$. Then $\gamma'$ is divided by $\TT$ into arc segments $\gamma_1$ and the arc segments of $g$, all of which are permissible. So $\gamma'$ is a permissible arc. Since $\gamma$ does not cross the interior of $\delta$, we have $\Int(\gamma',\gamma')=0=\Int(\gamma,\gamma')$. Moreover, we have $\chang{}{\gamma'}\leq\max\{j,m-(j-1)\}<m=\chang{}{\gamma}$ since $1<j<m$.
		
		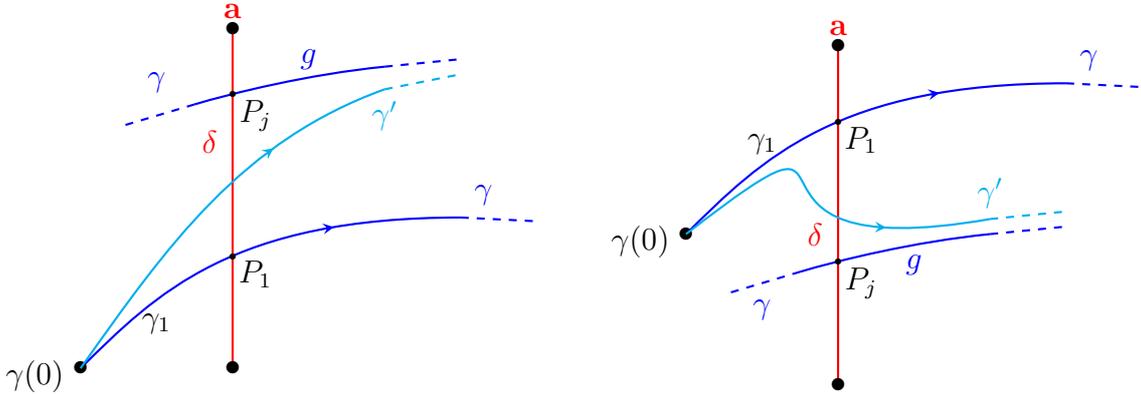
\begin{figure}[htpb]
			\begin{tikzpicture}
				\draw[thick,red] (5,0) .. controls (5,1) and (5,2) .. (5,4.5);
				\draw[thick,blue,->-=.7,>=stealth] (3,0)\nn .. controls (4,1) and (5,2) .. (8,2);
				\draw[thick,blue] (4.5,3.5) .. controls (5,3.65) and (6,3.9) .. (7,4);
				\draw[blue,thick,dashed] (7,4) to (8,4.1) (4.5,3.5) to (3.5,3.2) (8,2) to (9,1.95);
				\node at (2.4,-0.1){$\gamma(0)$};
				\draw (5,0)\nn (5,4.5)\nn;
				\draw[red](5,4.5)node[above]{$\ba$};
				\draw (5,3.6)node{$\bigcdot$};
				\node at (5.3,3.35){$P_j$};
				\draw[blue] (6,4.1)node{$g$};
				\draw (5,1.45) node{$\bigcdot$};
				\node at (5.3,1.25){$P_1$};
				\node at (4,.6){$\gamma_1$};
				\node[red] at (4.7,3){$\delta$};
				\draw[blue] (8.3,2.3)node{$\gamma$} (4,3.8)node{$\gamma$};
				\draw[cyan,thick,->-=.7,>=stealth] (3,0) ..controls +(55:2) and +(200:2) .. (7,3.7)node[below]{$\gamma'$};
				\draw[cyan, thick,dashed] (7,3.7)to(8,3.9);
			\end{tikzpicture}\qquad
			\begin{tikzpicture}
				\draw[thick,red] (5,0) .. controls (5,1) and (5,2) .. (5,4.5);
				\draw[thick,blue,->-=.7,>=stealth] (3,2)\nn .. controls (4,3) and (5,4) .. (8,4);
				\draw[thick,blue] (4.5,1.5) .. controls (5,1.65) and (6,1.9) .. (7,2);
				\draw[blue,thick,dashed] (7,2) to (8,2.1) (4.5,1.5) to (3.5,1.2) (8,4) to (9,3.95);
				\node at (2.4,1.9){$\gamma(0)$};
				\draw (5,0)\nn (5,4.5)\nn;
				\draw[red](5,4.5)node[above]{$\ba$};
				\draw (5,1.6)node{$\bigcdot$};
				\node at (5.3,1.35){$P_j$};
				\draw[blue] (6,1.6)node{$g$};
				\draw (5,3.45)node{$\bigcdot$};
				\node at (5.3,3.25){$P_1$};
				\node at (4,3.2){$\gamma_1$};
				\node[red] at (4.7,2){$\delta$};
				\draw[blue] (8.3,4.3)node{$\gamma$} (4,1)node{$\gamma$};
				\draw[cyan,thick,->-=.7,>=stealth] (3,2) ..controls +(38:3.5) and +(190:4) .. (7,2.2)node[above]{$\gamma'$};
				\draw[cyan, thick,dashed] (7,2.2)to(8,2.3);
			\end{tikzpicture}
			\caption{Case 1 in the proof of Lemma~\ref{l:sincere-reachable}}\label{f:subcase1.1}
		\end{figure}

		\noindent{\bf Case 2: $\ba_1=\ba_m$ does not contain $P_j$ for any $1<j<m$, and $\gamma_1$ and $\gamma_{m+1}$ are on the different sides of the arc $\ba_1=\ba_m$.}
		
		Let $\delta$ be the segment of $\ba_1=\ba_m$ from $P_1$ to $P_m$. See Figure~\ref{f:case3}. Take $\gamma'=\gamma_1\circ\delta\circ\gamma_{m+1}$. Then the arc segments of $\gamma$ divided by $\TT$ are $\gamma_1$ and $\gamma_{m+1}$, both of which are permissible. So $\gamma'$ is a permissible arc. Since $\gamma$ does not cross the interior of $\delta$, we have $\Int(\gamma',\gamma')=0=\Int(\gamma,\gamma')$. Moreover, we have $\chang{}{\gamma'}=1<m=\chang{}{\gamma}$.

		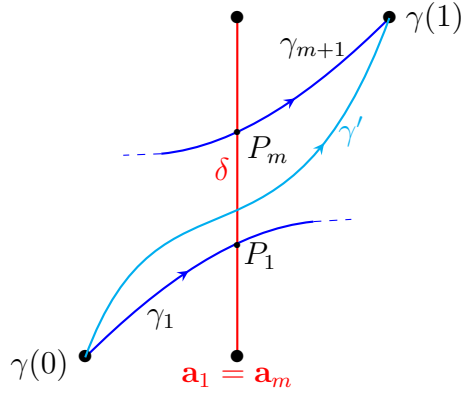
\begin{figure}[htpb]
			\begin{tikzpicture}
				\draw[thick,red] (5,0) .. controls (5,1) and (5,2) .. (5,4.5);
				\draw[thick,blue,->-=.5,>=stealth] (3,0)\nn .. controls (4,1) and (5,1.7) .. (6,1.8);
				\draw[thick,blue,-<-=.5,>=stealth] (7,4.5) .. controls (6,3.5) and (5, 2.8) .. (4,2.7);
				\node at (2.4,-0.1){$\gamma(0)$};
				\draw (7,4.5)\nn;
				\node at (7.6,4.5){$\gamma(1)$};
				\draw (5,0)\nn (5,4.5)\nn;
				\draw[red] (5,0)node[below]{$\ba_1=\ba_m$};
				\draw (5,2.95)node{$\bigcdot$};
				\node at (5.4,2.75){$P_m$};
				\draw (5,1.45)node{$\bigcdot$};
				\node at (5.3,1.35){$P_1$};
				\node[red] at (4.8,2.5){$\delta$};
				\node at (4,.5){$\gamma_1$};
				\node at (6,4.1){$\gamma_{m+1}$};
				\draw[cyan,thick,->-=.7,>=stealth] (3,0)..controls +(70:3) and +(-110:4)..(7,4.5);
				\draw[cyan] (6.5,3)node{$\gamma'$};
				\draw[blue,dashed] (6,1.8)to(6.5,1.83) (4,2.7)to(3.5,2.68);
			\end{tikzpicture}
			\caption{Case 2 in the proof of Lemma~\ref{l:sincere-reachable}}\label{f:case3}
		\end{figure}
		
		\noindent{\bf Case 3: both $\ba_1$ and $\ba_m$ do not contain $P_j$ for any $1<j<m$, and if $\ba_1=\ba_m$ then $\gamma_1$ and $\gamma_{m+1}$ are on the same side of $\ba_1$.}
		
		Reversing the orientation of $\gamma$ if necessary, we may assume that if $\ba_1=\ba_m$, then $P_m$ is in the interior of $\delta_1'$. So we have that $P_1$ is left most. By Lemma~\ref{l:concatenation}, this case can be divided into the following two subcases 3.1 and 3.2.

		\noindent{\bf Subcase 3.1: there exists $1<j<m+1$ such that $\delta_{j-1}'$, $\delta_j'$ and $\gamma_j$ are the three edges of $\triangle(\gamma_j)$.}

		Choose the smallest $j$ satisfying this condition. See Figure~\ref{f:subcase 2.1}. Since $P_1$ is left most, using Lemma~\ref{l:concatenation} repeatedly, we have that $P_1,\cdots, P_{j-1}$ are left most and the concatenation $\delta_{j-1}\circ\gamma_j$ is a permissible arc segment satisfying condition (P1). Since $\gamma$ has no self-intersections and the interior $\delta_{j-1}$ does not cross $\gamma$, the interior of $\delta_{j-1}\circ\gamma_j$ does not cross itself and $\gamma$. Take $\gamma'=(\delta_{j-1}\circ\gamma_j)\circ\gamma_{j+1}\circ\dots\circ\gamma_{m+1}$. Then $\gamma'$ is a permissible arc satisfying $\Int(\gamma',\gamma')=0=\Int(\gamma',\gamma)$. Moreover, we have $\chang{}{\gamma'}=m-(j-1)<m=\chang{}{\gamma}$ since $j>1$.
		
		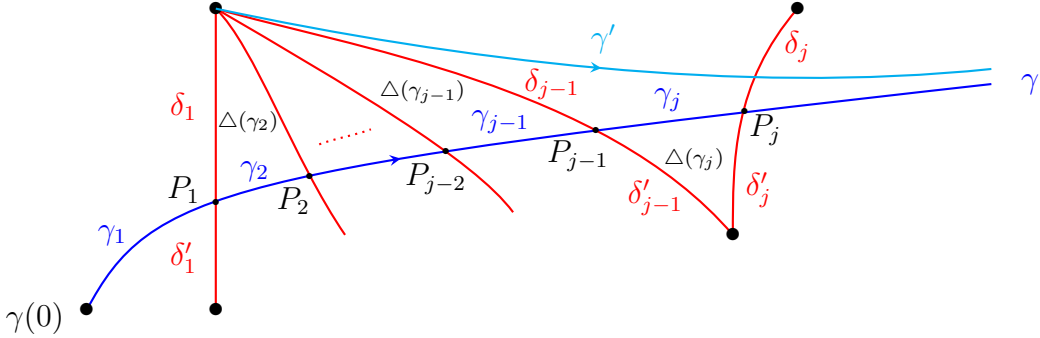
\begin{figure}[htpb]
			\begin{tikzpicture}[xscale=1.7]
				\draw[thick,red] (5,0) .. controls (5,1) and (5,2) .. (5,4);
				\draw[thick,blue,->-=.38,>=stealth] (4,0) .. controls (4.5,1.6) and (5,1.85) .. (11,3);
				\draw[thick,red] (5,4) .. controls (5.3,3.5) and (5.6,2) .. (6,1);
				\draw[thick,red] (5,4) .. controls (6,3) and (7,2) .. (7.3,1.3);
				\draw[thick,red] (5,4) .. controls (6,3.5) and (8,3) .. (9,1);
				\draw[thick,red] (9,1) .. controls (9,2) and (9,3) .. (9.5,4)\nn;
				\draw[thick,red,dotted] (5.8,2.2) .. controls (5.9,2.27) and (6.1,2.33) .. (6.2,2.4);
				\draw (9,1)\nn (4,0)\nn (5,0)\nn (5,4)\nn;
				\node[red] at (8.4,1.5){$\delta_{j-1}'$};
				\node at (3.6,-0.1){$\gamma(0)$};
				\draw (5,1.4)node{$\bigcdot$};
				\node at (4.75,1.6){$P_1$};
				\node[red] at (4.75,2.7){$\delta_1$};
				\node[red] at (4.75,0.7){$\delta_1'$};
				\node[blue] at (4.2,1){$\gamma_1$};
				\draw (5.725,1.75)node{$\bigcdot$};
				\node at (5.6,1.45){ $P_2$};
				\draw (6.78,2.07)node{$\bigcdot$};
				\node at (6.7,1.74){$P_{j-2}$};
				\draw (7.94,2.35)node{$\bigcdot$};
				\node at (7.8,2.05){$P_{j-1}$};
				\draw (9.09,2.6)node{$\bigcdot$};
				\node at (9.25,2.4){$P_{j}$};
				\node[blue] at (5.3,1.85){$\gamma_{2}$};
				\node[blue] at (7.2,2.5){$\gamma_{j-1}$};
				\node[blue] at (8.5,2.79){$\gamma_{j}$};
				\node[red] at (7.6,3){$\delta_{j-1}$};
				\draw[red] (9.2,1.6)node{$\delta'_j$} (9.5,3.5)node{$\delta_j$};
				\node[blue] at (11.3,3){$\gamma$};
				\draw (8.7,2)node{\tiny $\triangle(\gamma_j)$} (6.6,2.9)node{\tiny $\triangle(\gamma_{j-1})$} (5.25,2.5)node{\tiny $\triangle(\gamma_{2})$};
				\draw[cyan,thick,->-=.5,>=stealth] (5,4)..controls +(-20:2) and +(190:2) .. (11,3.2);
				\draw[cyan] (8,3.6)node{$\gamma'$};
			\end{tikzpicture}
			\caption{Subase 3.1}\label{f:subcase 2.1}
		\end{figure}
		
		\noindent{\bf Subcase 3.2: for any $1<j<m+1$, $\delta_{j-1}$, $\delta_j$ and $\gamma_j$ are the three edges of $\triangle(\gamma_j)$.}
		
		Then by Lemma~\ref{l:concatenation}~(b), all $P_j,1\leq j\leq m$, are left most. We denote by $O$ the common endpoint of $\delta_1,\ldots,\delta_m$. See Figure~\ref{f:subcase2.2}. This implies that for any $\ba\in\TT$, we have $\Int(\gamma,\ba)\leq 2$. 
		
		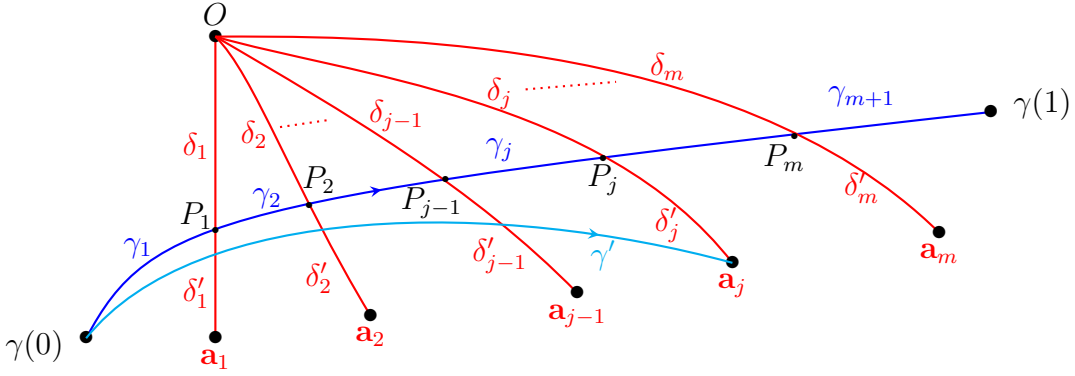
\begin{figure}[htpb]
			\begin{tikzpicture}[xscale=1.7]
				\draw[thick,red] (5,0)\nn .. controls (5,1) and (5,2) .. (5,4)\nn (5,0)node[below]{$\ba_1$};
				\draw[thick,blue,->-=.36,>=stealth] (4,0)\nn .. controls (4.5,1.6) and (5,1.85) .. (11,3)\nn;
				\draw[thick,red] (5,4) .. controls (5.3,3.5) and (5.6,2) .. (6.2,0.3)\nn (6.2,.3)node[below]{$\ba_2$};
				\draw[thick,red] (5,4) .. controls (6,3) and (7,2) .. (7.8,0.6)\nn (7.8,0.6)node[below]{$\ba_{j-1}$};
				\draw[thick,red] (5,4) .. controls (6,3.5) and (8,3.2) .. (9,1)\nn (9,1)node[below]{$\ba_j$};
				\draw[thick,red,dotted] (5.5,2.8) to (5.9,2.9);
				\draw[thick,red,dotted] (7.4,3.3) to (8.1,3.4);
				\draw[thick,red] (5,4) .. controls (6,4) and (9,4.1) .. (10.6,1.4)\nn (10.6,1.4)node[below]{$\ba_m$};
				\node at (3.6,-0.1){$\gamma(0)$};
				\node at (11.4,3.1){$\gamma(1)$};
				\node at (5,4.3){$O$};
				\draw (5,1.4)node{$\bigcdot$};
				\node at (4.85,1.63){$P_1$};
				\node[red] at (4.85,2.55){$\delta_1$};
				\node[red] at (5.3,2.7){$\delta_2$};
				\node[red] at (6.4,2.95){$\delta_{j-1}$};
				\node[red] at (7.2,3.25){$\delta_{j}$};
				\node[red] at (8.5,3.6){$\delta_{m}$};
				\node[red] at (4.85,0.6){$\delta_1'$};
				\node[red] at (5.8,0.8){$\delta_2'$};
				\node[red] at (7.2,1.1){$\delta_{j-1}'$};
				\node[red] at (8.5,1.5){$\delta_{j}'$};
				\node[red] at (10,2){$\delta_{m}'$};
				\node[blue] at (4.4,1.2){$\gamma_1$};
				\draw (5.725,1.73)node{$\bigcdot$};
				\node at (5.8,2.1){$P_2$};
				\draw (6.78,2.07)node{$\bigcdot$};
				\node at (6.67,1.8){$P_{j-1}$};
				\draw (8,2.35)node{$\bigcdot$};
				\node at (8,2.1){$P_{j}$};
				\draw (9.48,2.65)node{$\bigcdot$};
				\node at (9.4,2.4){$P_{m}$};
				\node[blue] at (5.4,1.85){$\gamma_{2}$};
				\node[blue] at (7.2,2.5){$\gamma_{j}$};
				\node[blue] at (10,3.2){$\gamma_{m+1}$};
				\draw[cyan,thick,->-=.8,>=stealth] (4,0) .. controls +(65:2) and +(155:2) .. (9,1);
				\draw[cyan,thick] (8,1.1)node{$\gamma'$};
			\end{tikzpicture}
			\caption{Subcase 3.2}\label{f:subcase2.2}
		\end{figure}
		
		\noindent{\bf Subcase 3.2.1: there exists $2\leq j\leq m$ such that $\ba_j\neq\ba_i$ for any $1\leq i\leq m$ and $i\neq j$.}
		
		Take $\gamma'=\gamma_1\circ\gamma_2\circ\dots\circ\gamma_{j-1}\circ(\gamma_j\circ{\delta'_j}^{-1}),$ see Figure~\ref{f:subcase2.2}. Applying Lemma~\ref{l:concatenation}~(b) to the inverse of $\gamma$, we have that $\gamma_j\circ{\delta_j'}^{-1}$ is a permissible arc segment satisfying condition (P1).
		So $\gamma'$ is a permissible arc. Since the interior of $\delta_j'$ does not cross $\gamma$, a similar discussion as in Subcase 3.1 shows that $\Int(\gamma',\gamma')=0=\Int(\gamma',\gamma)$. Moreover, we have $\chang{}{\gamma'}=j-1<m=\chang{}{\gamma}$.

		\noindent{\bf Subcase 3.2.2: for any $2\leq j\leq m$, there is $1\leq i\leq m$ and $i\neq j$ such that $\ba_i=\ba_j$.}
		
		Since $\ba_m$ does not contain $P_2,\cdots,P_{m-1}$, we have $\ba_m=\ba_1$. So for any $1\leq j\leq m$, there is $i\neq j$ such that $\ba_i=\ba_j$. Recall that $0<\Int(\gamma,\ba)\leq 2$ holds for any $\ba\in\TT$ in Subcase 3.2. So in Subcase 3.2.2, we have $\Int(\gamma,\ba)= 2$ for any $\ba\in\TT$, and $m=2|\TT|$. It follows that all arcs in $\TT$ are loops sharing the same endpoint, say $O$, and $\gamma$ cuts out each angle formed by two arcs in $\TT$ exactly once. Recall that $\gamma_1$ and $\gamma_{m+1}$ are on the same side of $\ba_m=\ba_1$. It follows that the tile containing $\gamma_1$ is the same as the tile containing $\gamma_{m+1}$, none of whose angles is formed by two (not necessarily distinct) arcs in $\TT$. So this tile is of type (IV), i.e. a $3$-gon whose three edges are $\ba_1=\ba_m$ and two boundary components with $O$ a common endpoint and without any unmarked boundary component in its interior. Thus, the boundary component of $\SS$ containing $O$ has two marked points, say $P$ the other marked point. Then we have $\gamma(0)=\gamma(1)=P$. See Figure~\ref{f:subcase:3.2.2}.
		
		\begin{figure}[htpb]
			\begin{tikzpicture}[xscale=2]
				\draw[thick,fill=black!20] (0,1)ellipse(.5 and 1);
				\draw[thick,red] (0,0) .. controls +(-10:2) and +(0:2) .. (0,4) (0,0) .. controls +(190:2) and +(180:2) .. (0,4)node[below]{$\ba_1=\ba_m$};
				\draw[red,thick] (0,0) to (-3,-3)node[below]{$\ba_2$} (0,0) to (3,-3)node[left]{$\ba_{m-1}$};
				\draw[blue,thick] (0,2) .. controls +(20:3) and +(0:3) .. (0,-2);
				\draw[blue,thick,->-=.5,>=stealth] (0,2) .. controls +(160:3) and +(180:3) .. (0,-2);
				\draw[blue] (.8,2.4)node{$\gamma_{m+1}$} (-.8,2.4)node{$\gamma_1$} (-2.5,0)node{$\gamma_2$} (2.3,0)node{$\gamma_m$};
				\draw[red,thick] (0,0)to(-1,-4)node[right]{$\ba_j$};
				\draw[red,dotted,thick] (-125:1)to(-110:1) (-70:1)to(-55:1);
				\draw (-1.5,0)node{$\triangle(\gamma_2)$} (1.5,0)node{$\triangle(\gamma_m)$} (-1,1.4)node{$\triangle(\gamma_1)$} (1,1.4)node{$\triangle(\gamma_{m+1})$};
				\draw (0,2)node[above]{$P$} (0,0)node[below]{$O$};
				\draw (-1.63,2.17)node{$P_1$} (1.63,2.17)node{$P_m$} (-1.5,-1.9)node{$P_2$} (1.5,-1)node{$P_{m-1}$} (-.4,-2.3)node{$P_{j}$} (.38,-2.3)node{$P_{j+1}$};
				\draw[blue] (0,-1.75)node{$\gamma_{j+1}$};
				\draw[red,thick,dashed] (3,-3)..controls+(-45:1.5) and +(-104:1.5) .. (-1,-4);
				\draw[red] (2.3,-4.3)node{$\delta$};
				\draw[cyan,thick,->-=.5,>=stealth] (0,2) .. controls +(50:4) and +(45:2) .. (1.8,-1.8) ..controls +(-135:1) and +(-14:1) .. (-.65,-2.6) ..controls +(166:.3) and +(140:.6).. (-1,-4.8);
				\draw[cyan,dashed,thick]  (-1,-4.8) .. controls +(-30:.8) and +(-140:3) .. (3.4,-3.8);
				\draw[cyan,thick] (3.4,-3.8) ..controls +(40:.8) and +(-45:.6).. (2.3,-2) ..controls +(135:.1) and +(-115:.3).. (2.3,-1.4) ..controls +(65:.6) and +(-90:1).. (2.7,1)..controls +(90:1) and +(55:3.5)..(0,2);
				\draw[cyan] (2.5,3)node{$\gamma'$};
				\draw[red,thick] (0,0)to(1,-4)node[right]{$\ba_{j+1}$};
				\draw (-1.485,1.855)node{$\bigcdot$} (1.485,1.855)node{$\bigcdot$} (-1.46,-1.49)node{$\bigcdot$} (1.46,-1.49)node{$\bigcdot$} (-.485,-2)node{$\bigcdot$} (.485,-2)node{$\bigcdot$};
				\draw (0,2)\nn (0,0)\nn;
			\end{tikzpicture}
			\caption{Subcase 3.2.2}\label{f:subcase:3.2.2}
		\end{figure}
		
		Since $\Int(\gamma,\ba)= 2$ for any $\ba\in\TT$, there is $1<j<m-1$ such that $\ba_j=\ba_{m-1}$. Let $\delta$ be the segment of $\ba_j=\ba_{m-1}$ from $P_j$ to $P_{m-1}$. The concatenation $\delta\circ\gamma_m$ is isotopic to $\gamma_m$ relative to the interiors of arcs in $\TT$. So it is permissible. Take $\gamma'=\gamma_{m+1}^{-1}\circ\gamma_m^{-1}\circ\cdots\circ\gamma_{j+1}^{-1}\circ(\delta\circ\gamma_m)\circ\gamma_{m+1}$. Then $\gamma'$ is a permissible curve with $\chang{}{\gamma'}=m-j+2$. Since $P_j$ and $P_{m-1}$ are the only intersections between $\gamma$ and $\ba_j=\ba_{m-1}$, the interior of $\delta$ does not cross $\gamma$. So we have $\Int(\gamma',\gamma')=0$ and $\Int(\gamma',\gamma)=0$.
		
		If $j>2$, then $\chang{}{\gamma'}=m-j+2<m=\chang{}{\gamma}$ and we are done. If $j=2$, then $\chang{}{\gamma'}=m=\chang{}{\gamma}$. We need to find another permissible arc $\gamma''$ without self-intersections such that $\Int(\gamma'',\gamma')=0$ and $\chang{}{\gamma''}<\chang{}{\gamma'}$. Note that in case $j=2$, $\gamma_2$ and $\gamma_m$ are in the same tile, where both of the angles incident to $\ba_2=\ba_{m-1}$ have $\ba_1=\ba_m$ as an edge. So this tile is of type II, whose edges are $\ba_1=\ba_m$ and $\ba_2=\ba_{m-1}$ and which contains an unmarked boundary component in its interior. See Figure~\ref{f:subcase:3.2.2.1}. Let $\gamma''$ be the concatenation of $\gamma'$ with the boundary segment edge of $\triangle(\gamma_1)$ from $P$ to $O$. By definition, the last arc segment of $\gamma''$ is permissible in a tile of type (II). So $\gamma''$ is permissible. Moreover, we have $\Int(\gamma'',\gamma'')=0$, $\Int(\gamma'',\gamma')=0$ and $\chang{}{\gamma''}=m-1<m=\chang{}{\gamma'}$ as required.
		
		\begin{figure}[htpb]
			\begin{tikzpicture}[xscale=2]
				\draw[thick,fill=black!20] (0,1)ellipse(.5 and 1);
				\draw[thick,fill=black!20] (0,4.6)ellipse(.15 and .3);
				\draw (0,2)\nn (0,0)\nn;
				\draw[thick,red] (0,0) .. controls +(-10:2) and +(0:2) .. (0,4) (0,0) .. controls +(190:2) and +(180:2) .. (0,4)node[below]{$\ba_1=\ba_m$};
				\draw[thick,red] (0,0) .. controls +(-20:3) and +(0:3) .. (0,6) (0,0) .. controls +(200:3) and +(180:3) .. (0,6)node[below]{$\ba_2=\ba_{m-1}$};
				\draw[red,thick] (0,0) to (-2.25,-2.25)node[below]{$\ba_3$} (0,0) to (2.25,-2.25)node[below]{$\ba_{m-2}$};
				\draw[blue,thick] (0,2) .. controls +(20:3) and +(0:3) .. (0,-2);
				\draw[blue,thick,->-=.5,>=stealth] (0,2) .. controls +(160:3) and +(180:3) .. (0,-2);
				\draw[blue] (1,2.35)node{$\gamma_{m+1}$} (-.8,2.4)node{$\gamma_1$};
				\draw[red,dotted,thick] (-125:1)to(-55:1);
				\draw (-1,1.4)node{$\triangle(\gamma_1)$} (1,1.4)node{$\triangle(\gamma_{m+1})$};
				\draw (0,2)node[above]{$P$} (0,0)node[below]{$O$};
				\draw (-1.63,2.17)node{$P_1$} (1.63,2.17)node{$P_m$} (-1.5,-1.9)node{$P_2$} (1.5,-1)node{$P_{m-1}$};
				\draw[cyan,thick,->-=.42,>=stealth] (0,2) .. controls +(55:4.5) and +(45:2) .. (2.05,-2.05) ..controls +(-135:1) and +(0:.5) .. (0,-2.95) ..controls +(180:.5) and +(-45:1).. (-2.05,-2.05)..controls +(135:2) and +(180:2.4)..(0,5.4) ..controls +(0:1.5) and +(60:1) ..(0,2);
				\draw[green,thick,->-=.35,>=stealth] (0,2) .. controls +(50:4) and +(45:2) .. (1.9,-1.9) ..controls +(-135:1) and +(0:.5) .. (0,-2.7) ..controls +(180:.5) and +(-45:1).. (-1.9,-1.9)..controls +(135:2) and +(180:2)..(0,5.1) ..controls +(0:.5) and +(0:.5) ..(0,4.1)..controls +(180:3) and +(200:2).. (0,0);
				\draw[cyan] (0,-3.2)node{$\gamma'$};
				\draw[green] (0,-2.45)node{$\gamma''$};
				\draw (-1.485,1.855)node{$\bigcdot$} (1.485,1.855)node{$\bigcdot$} (-1.98,1.25)node{$\bigcdot$} (1.98,1.25)node{$\bigcdot$} (-1.46,-1.49)node{$\bigcdot$} (1.46,-1.49)node{$\bigcdot$} (-2,1.25)node[left]{$P_2$} (1.95,1.25)node[right]{$P_{m-1}$};
			\end{tikzpicture}
			\caption{$\ba_2=\ba_{m-1}$ in Subcase 3.2.2}\label{f:subcase:3.2.2.1}
		\end{figure}
	\end{proof}

	\begin{theorem}\label{t:t-reachability}
		Let $A$ be a finite-dimensional gentle algebra over $k$ such that $|A|>1$. Then $A$  has the $\tau$-reachable property.
	\end{theorem}
	\begin{proof}
		It suffices to show that every indecomposable $\tau$-rigid pair $(M,0)$ is $\tau$-reachable from $(0,P)$ for some indecomposable projective $A$-module $P$. This is clear whenever $M$ is non-sincere. Now assume that $M$ is sincere, applying Lemma~\ref{l:sincere-reachable} repeatedly, we know that $M$ is $\tau$-reachable from a non-sincere indecomposable $\tau$-rigid $A$-module $N$. Consequently, $(M,0)$ is $\tau$-reachable from $(0,P)$ for some indecomposable projective $A$-module $P$.
	\end{proof}

	\section{Proofs of the main results}\label{s:proofs-main-results}
	
	Let $A$ be a finite-dimensional $k$-algebra. Let $(U,R)$ be a basic $\tau$-rigid pair. We refer to Appendix~\ref{s:reduction} for the reduction theory of $(U,R)$. Denote by $(T_{(U,R)},R)$ the Bongartz completion of $(U,R)$. Let $A_{(U,R)}:=\End_A(T_{(U,R)})/\langle e_U\rangle$ be the factor algebra of $\End_A(T_{(U,R)})$ by the ideal generated by $e_U$, where $e_U$ is the idempotent of $\End_A(T_{(U,R)})$ associated to the direct summand $U$.
	
	\begin{definition}\label{d:totally-tau-reachable}
		We call $A$ has the totally $\tau$-reachable property provided that for any basic $\tau$-rigid pair $(U,R)$, the algebra $A_{(U,R)}$ is $\tau$-reachable when $|A_{(U,R)}|>1$.
	\end{definition}
	
	By definition, if $A$ is totally $\tau$-reachable and $|A|>1$, then $A$ is $\tau$-reachable; if $|A|=1$, then $A$ is always totally $\tau$-reachable.

	\begin{lemma}\label{l:transitivity}
		If $A$ is totally $\tau$-reachable, then for any basic $\tau$-rigid pair $(U,R)$, the algebra $A_{(U,R)}$ is also totally $\tau$-reachable.
	\end{lemma}
	
	\begin{proof}
		For any basic $\tau$-rigid pair $(V,Q)$ in $\mod A_{(U,R)}$, by Corollary~\ref{c:useful}, there is a basic $\tau$-rigid pair $(V',Q')$ in $\mod A$ such that $(A_{(U,R)})_{(V,Q)}\cong A_{(V',Q')}$. Since $A$ is totally $\tau$-reachable, $A_{(V',Q')}$ is $\tau$-reachable if $|A_{(V',Q')}|>1$. This implies that $A_{(U,R)}$ is totally $\tau$-reachable.
	\end{proof}
	
	\begin{proposition}\label{p:totally-tau-reachable}
		Let $A$ be a finite-dimensional $k$-algebra. If $A$ is totally $\tau$-reachable, then the support $\tau$-tilting graph $\mh(\opname{s\tau-tilt} A)$ of $A$ is connected.
	\end{proposition}
	
	\begin{proof}
		Use induction on $|A|$. If $|A|=1$, by Remark~\ref{rmk:rank1}, the support $\tau$-tilting graph of $A$ is connected. Assume that the theorem holds for the case $|A|<n$. Now suppose $|A|=n$. Let $(M,P)$ and $(N,Q)$ be two basic support $\tau$-tilting pairs of $A$-modules. We separate the remaining proof by considering whether $(M,P)$ and $(N,Q)$ share a non-zero direct summand.
		
		\noindent{\bf Case 1:} $(M,P)$ and $(N,Q)$ have a common direct summand $(U,R)$. By Lemma~\ref{l:transitivity}, $A_{(U,R)}$ is  totally $\tau$-reachable with $|A_{(U,R)}|<|A|=n$. By the inductive hypothesis, the support $\tau$-tilting graph $\mh(\opname{s\tau-tilt} A_{(U,R)})$ of $A_{(U,R)}$ is connected. By Corollary~\ref{c:key}, the graph $\mh(\opname{s\tau-tilt} A_{(U,R)})$ is isomorphic to the full subgraph of $\mh(\opname{s\tau-tilt}A)$ whose vertices are the support $\tau$-tilting pairs containing $(U,R)$ as a direct summand. So $(M,P)$ and $(N,Q)$ are in this subgraph. It follows that they are mutation-reachable to each other.

		\noindent{\bf Case 2:} $(M,P)$ and $(N,Q)$ do not have any common nonzero direct summand. Let $(M_1, P_1)$ be an indecomposable direct summand of $(M,P)$ and $(N_1, Q_1)$ an indecomposable direct summand of $(N,Q)$. Since $A$ is $\tau$-reachable, $(M_1, P_1)$ is $\tau$-reachable from $(N_1, Q_1)$. By definition, there exists a sequence of $\tau$-rigid pairs $(L_1,R_1),\ldots, (L_t,R_t)$ such that $(L_i\oplus L_{i+1}, R_i\oplus R_{i+1})$, $0\leq i\leq t$, are $\tau$-rigid pairs in $\mod A$, where $(L_0,R_0)=(M_1,P_1)$ and $(L_{t+1},R_{t+1})=(N_1,Q_1)$. According to \cite[Theorem 2.10]{AIR}, each $\tau$-rigid pair $(L_i\oplus L_{i+1}, R_i\oplus R_{i+1})$ can be completed into a support $\tau$-tilting pair $(\tilde{L}_{i}, \tilde{R}_{i})$. Denote $(\tilde{L}_{-1},\tilde{R}_{-1})=(M,P)$ and $(\tilde{L}_{t+1},\tilde{R}_{t+1})=(N,Q)$. Then
		$(\tilde{L}_i,\tilde{R}_i)$ and $(\tilde{L}_{i+1}, \tilde{R}_{i+1})$ have a common direct summand $(L_{i+1},R_{i+1})$ for $-1\leq i\leq t$. We conclude that $(M,P)$ is mutation-reachable by $(N,Q)$ by using the result in Case 1 repeatedly. This completes the proof.
	\end{proof}

	\begin{corollary}\label{c:reachable-in-face}
		Let $A$ be a finite-dimensional $k$-algebra. If $A$ is totally $\tau$-reachable, then $A$ has the reachable-in-face property.
	\end{corollary}
	
	\begin{proof}
		By Lemma~\ref{l:transitivity}, for any basic $\tau$-rigid pair $(U,R)$ in $\mod A$, the algebra $A_{(U,R)}$ is totally $\tau$-reachable. Then by Proposition~\ref{p:totally-tau-reachable}, the support $\tau$-tilting graph of $A_{(U,R)}$ is connected. This graph is isomorphic to the face $\mathcal{F}_{(U,R)}$. So any face of $\mh(\opname{s\tau-tilt} A)$ is connected. In particular, $A$ has the reachable-in-face property.
	\end{proof}

	\begin{theorem}\label{t:totally-equivalent-to-connected}
		Let $A$ be a finite-dimensional $k$-algebra. Then $A$ is totally $\tau$-reachable if and only if $A$ has the reachable-in-face property and the support $\tau$-tilting graph $\mh(\opname{s\tau-tilt} A)$ is connected.
	\end{theorem}
	\begin{proof}
		The ``only if" part follows from Proposition \ref{p:totally-tau-reachable} and Corollary \ref{c:reachable-in-face}. Let us prove the ``if" part. Assume that $A$ has the reachable-in-face property and the support $\tau$-tilting graph $\mh(\opname{s\tau-tilt} A)$ is connected. Consequently, each face of $\mh(\opname{s\tau-tilt} A)$ is connected.
		Let $(U,R)$  be a basic $\tau$-rigid pair of $A$ such that $|A_{(U,R)}|>1$, we have to show that $A_{(U,R)}$ is $\tau$-reachable. According to Corollary \ref{c:key}, there is a bijection $E_{(U,R)}:\opname{s\tau-tilt-pair}_{(U,R)}A\to\opname{s\tau-tilt-pair} A_{(U,R)}$ which commutes with the mutation. If follows that the support $\tau$-tilting graph $\mh(\opname{s\tau-tilt} A_{(U,R)})$ of $A_{(U,R)}$ is isomorphic to the face $\mathcal{F}_{(U,R)}$. In particular, $\mh(\opname{s\tau-tilt} A_{(U,R)})$ of $A_{(U,R)}$ is connected. Now the result follows from Proposition \ref{p:connect-reachable}.
	\end{proof}
	
	\begin{proposition}\label{p:gentle-are-totally-tau-reachable}
		Any finite-dimensional gentle $k$-algebra  is totally $\tau$-reachable.
	\end{proposition}
	
	\begin{proof}
		Let $A$ be a finite-dimensional gentle algebra over $k$. For any basic $\tau$-rigid pair $(U,R)$ in $\mod A$, by Lemma~\ref{l:endomorphism algebra is gentle}, 
		the endomorphism algebra $\End_A(T_{(U,R)})$ is gentle. By Lemma~\ref{l:factor-gentle-algebra}, the factor algebra $A_{(U,R)}$ of $\End_A(T_{(U,R)})$ by the ideal generated by an idempotent is also gentle. Hence by Theorem~\ref{t:t-reachability}, $A_{(U,R)}$ is $\tau$-reachable if $|A_{(U,R)}|>1$. Thus, $A$ is totally $\tau$-reachable.
	\end{proof}
	
	Combining Proposition~\ref{p:gentle-are-totally-tau-reachable} and Theorem~\ref{t:totally-equivalent-to-connected}, we get Theorem~\ref{t:main-thm-1} and Theorem~\ref{t:main-thm-2}.

	\appendix
	
	\newcommand{\tor}{\mathcal{T}}
	\newcommand{\tof}{\mathcal{F}}

	\section{Reduction}\label{s:reduction}
	
	\newcommand{\W}{\mathcal{W}}

	Let $A$ be a finite-dimensional $k$-algebra and $\mod A$ the category of finitely generated right $A$-modules. For any $M\in\mod A$, denote by $\opname{Fac} M$ the full subcategory of $\mod A$ consisting of all factors modules of direct sums of copies of $M$; denote $$M^\perp=\{N\in \mod A~|~ \Hom_A(M,N)=0\},\ ^\perp M=\{N\in \mod A~|~ \Hom_A(N,M)=0\}.$$
	
	Denote by $\opname{\tau-rigid-pair}A$ the set of isoclasses of basic $\tau$-rigid pairs in $\mod A$ and by $\opname{s\tau-tilt-pair}A$ the set of isoclasses of basic support $\tau$-tilting pairs in $\mod A$. There is a partial order on $\opname{s\tau-tilt-pair}A$ that $(M,P)\leq (N,Q)$ if and only if $\Fac M\subseteq\Fac N$. There is a bijection from $\opname{s\tau-tilt-pair}A$ to the set of functorially finite torsion pairs in $\mod A$, sending $(M,P)$ to $(\tor_{M},\tof_M):=(\Fac M, M^\perp)$, see \cite[Theorem~2.7]{AIR}.

	Let $(U,R)$ be a basic $\tau$-rigid pair. Denote by $\opname{\tau-rigid-pair}_{(U,R)}A$ the set of isoclasses of $\tau$-rigid pairs containing $(U,R)$ as a direct summand, and by $\opname{s\tau-tilt-pair}_{(U,R)}A$ the set of isoclasses of support $\tau$-tilting pairs containing $(U,R)$ as a direct summand.

	\begin{lemma}[{\cite[Theorem 2.10]{AIR} and \cite[Theorem~4.4]{DIRRT}}]\label{l:AD}
		Let $(U,R)$ be a basic $\tau$-rigid pair in $\mod A$. There are $(S_{(U,R)},L_{(U,R)})\leq (T_{(U,R)},R)\in \opname{s\tau-tilt-pair}_{(U,R)}A$ such that $\opname{s\tau-tilt-pair}_{(U,R)}A$ is the interval of $\opname{s\tau-tilt-pair}A$ between $(S_{(U,R)},L_{(U,R)})$ and $(T_{(U,R)},R)$. Moreover, the corresponding torsion pairs are $$(\tor_{T_{(U,R)}},\tof_{T_{(U,R)}})=({}^{\perp}(\tau U)\cap R^\perp,T_{(U,R)}^\perp)\text{ and }(\tor_{S_{(U,R)}},\tof_{S_{(U,R)}})=(\Fac U,U^\perp).$$
	\end{lemma}
	
	The pair $(T_{(U,R)},R)$ (resp. $(S_{(U,R)},L_{(U,R)})$) is called the {\it Bongartz complement} (resp. {\it co-Bongartz complement}) of $(U,R)$. Let $\End_A(T_{(U,R)})$ be the endomorphism algebra of $T_{(U,R)}$ and $e_U$ the idempotent of $\End_A(T_{(U,R)})$ associated to the direct summand $U$. We denote by $A_{(U,R)}:=\End_A(T_{(U,R)})/\langle e_U\rangle$ the factor algebra of $\End_A(T_{(U,R)})$ by the ideal generated by $e_U$. We have $|A_{(U,R)}|=|A|-|U|-|R|$, which is called the {\it co-rank} of $(U,R)$, denoted by $\opname{co-rank}(U,R)$. Denote by
	$$\W(U,R)=\tor_{T_{(U,R)}}\cap \tof_{S_{(U,R)}}={}^{\perp}\tau U\cap R^\perp \cap U^{\perp}.$$

	\begin{lemma}[{\cite[Theorem~1.4]{J} (for $R=0$) and \cite[Theorem~4.12~(a,b)]{DIRRT}}]\label{l:wide}
		The subcategory $\W(U,R)$ is a wide subcategory of $\mod A$, and the functor $$\Hom_A(T_{(U,R)},-):\mod A\to \mod A_{(U,R)}$$ restricts to an equivalence
		$$F_{(U,R)}\colon\mathcal{W}(U,R)\to\mod A_{(U,R)}.$$
	\end{lemma}

	
	A pair $(M,P)$ of objects in $\mathcal{W}(U,R)$ is called a {\it $\tau$-rigid pair} if $(F_{(U,R)}(M),F_{(U,R)}(P))$ is a $\tau$-rigid pair in $\mod A_{(U,R)}$. Denote by $\opname{\tau-rigid-pair}\W(U,R)$ the set of isoclasses of basic $\tau$-rigid pairs in $\W(U,R)$. For any $(V,Q)\in\opname{\tau-rigid-pair}\W(U,R)$, denote by $\opname{\tau-rigid-pair}_{(V,Q)}\W(U,R)$ the subset of $\opname{\tau-rigid-pair}\W(U,R)$ consisting of the pairs containing $(V,Q)$ as a direct summand.
	
	For any $M\in\W(U,R)$, we denote by $\Fac_{\W(U,R)}M$ the full subcategory of $\W(U,R)$ consisting of objects $N$ such that there is an epimorphism in $\W(U,R)$ from a direct sum of copies of $M$ to $N$. The following result is essentially from \cite{BM}.
	
	\begin{proposition}\label{l:redBM}
		Let $(U,R)$ be a basic rigid pair in $\mod A$. Then there is a bijection
		$$\begin{array}{rccc}
			\mathcal{E}_{(U,R)}\colon&\opname{\tau-rigid-pair}_{(U,R)}A&\to&\opname{\tau-rigid-pair}\mathcal{W}(U,R)\\
			&(M,P)&\to&(\mathcal{E}'_{(U,R)}(M,P),\mathcal{E}''_{(U,R)}(M,P))\end{array}$$
		which commutes with direct sums and such that
		\begin{enumerate}
			\item for any $(M,P)\in\opname{\tau-rigid-pair}_{(U,R)}A$, we have $$\opname{co-rank}(M,P)=\opname{co-rank}(\mathcal{E}'_{(U,R)}(M,P),\mathcal{E}''_{(U,R)}(M,P)),$$
			\item for any $(M,P)\in \opname{s\tau-tilt-pair}_{(U,R)}A$, we have $$\Fac M\cap\W(U,R)=\Fac_{\W(U,R)}(\mathcal{E}'_{(U,R)}(M,P)).$$
		\end{enumerate}
	\end{proposition}
	
	\begin{proof}
		Applying \cite[Proposition~5.11]{BM} to $(0,R)$ in $\mod A$, we have that there is a bijection
		$$\begin{array}{rccc}
			\mathcal{E}_{(0,R)}\colon&\opname{\tau-rigid-pair}_{(0,R)}A&\to&\opname{\tau-rigid-pair}\mathcal{W}(0,R)\\
			&(M,P)&\to&(M,\mathcal{E}''_{(0,R)}(M,P))\end{array}$$
		which commutes with direct sums and such that for any $(M,P)\in\opname{\tau-rigid-pair}_{(0,R)}A$, we have $|P|=|\mathcal{E}''_{(0,R)}(M,P)|+|R|$. In particular, $(U,0)=\mathcal{E}_{(0,R)}((U,R))$ and $(T_{(U,R)},0)=\mathcal{E}_{(0,R)}((T_{(U,R)},R))$ are basic $\tau$-rigid pairs in $\mathcal{W}(0,R)$. Since the map $\mathcal{E}_{(0,R)}$ commutes with direct sums, it restricts to a bijection
		$$\begin{array}{ccc}
			\opname{\tau-rigid-pair}_{(U,R)}A&\to&\opname{\tau-rigid-pair}_{(U,0)}\mathcal{W}(0,R).
		\end{array}$$
		Note that $\W(0,R)=R^\perp$ is a Serre subcategory of $\mod A$, i.e., a full subcategory of $\mod A$ closed under extensions, factor modules and submodules. So for any $N\in\W(0,R)$, we have $\Fac N=\Fac_{\W(0,R)}N$. It follows that $(T_{(U,R)},0)$ is the Bongartz completion of $(U,0)$ in $\W(0,R)$. 
		
		Let $\W_{\W(0,R)}(U,0)$ be the wide subcategory of $\W(0,R)$ induced by the $\tau$-rigid pair $(U,0)$. Since $\Fac T_{(U,R)}={}^\perp \tau U\cap R^\perp\subseteq\W(0,R)$ and $\Fac T_{(U,R)}=\Fac_{\W(0,R)}T_{(U,R)}$, we have $$\W_{\W(0,R)}(U,0)=\Fac_{\W(0,R)} T_{(U,R)}\cap U^\perp\cap\W(0,R)=\Fac T_{(U,R)}\cap U^\perp=\W(U,R).$$ Applying \cite[Proposition~5.8]{BM} to $(U,0)$ in $\W(0,R)$, we have that there is a bijection
		$$\begin{array}{rccc}
			\mathcal{E}_{(U,0)}\colon&\opname{\tau-rigid-pair}_{(U,0)}\W(0,R)&\to&\opname{\tau-rigid-pair}\W_{\W(0,R)}(U,0)=\opname{\tau-rigid-pair}\mathcal{W}(U,R)\\
			&(M,P)&\mapsto&(\mathcal{E}'_{(U,0)}(M,P),\mathcal{E}''_{(U,0)}(M,P))
		\end{array}	$$
		which commutes with direct sums and such that for any $(M,P)\in\opname{\tau-rigid-pair}_{(U,0)}\W(0,R)$, we have $|M|+|P|=|\mathcal{E}'_{(U,0)}(M,P)|+|\mathcal{E}''_{(U,0)}(M,P)|+|U|$ and the map $\mathcal{E}'_{(U,0)}$ is the same as the map $f$ in \cite[Theorem~3.15]{J}. So for any $(M,P)\in \opname{s\tau-rigid-pair}_{(U,0)}\W(0,R)$ with $\opname{co-rank}(M,P)=0$, we have $\Fac M\cap\W(U,R)=\Fac_{\W(U,R)}\mathcal{E}'_{(U,0)}(M,P)$.
		
		Composing the two bijections $\mathcal{E}_{(0,R)}$ and $\mathcal{E}_{(U,0)}$, we get the required bijection.
	\end{proof}
	
	Combing the bijection in Proposition~\ref{l:redBM} with the equivalence in Lemma~\ref{l:wide}, we have the following reduction technique.

	\begin{corollary}\label{c:key}
		Let $(U,R)$ be a basic rigid pair in $\mod A$. There is a bijection
		$$\begin{array}{rccc}
			E_{(U,R)}\colon&\opname{\tau-rigid-pair}_{(U,R)}A&\to&\opname{\tau-rigid-pair} A_{(U,R)}\end{array}$$
		which commutes with direct sums and restricts to a bijection
		$$\begin{array}{rccc}
			E_{(U,R)}\colon&\opname{s\tau-tilt-pair}_{(U,R)}A&\to&\opname{s\tau-tilt-pair} A_{(U,R)}\end{array}$$
		which preserves the order and commutes with the mutation, and such that for any $(M,P)\in \opname{s\tau-tilt-pair}_{(U,R)}A$, we have 
		$$F_{(U,R)}\left(\tor_M\cap\W(U,R)\right)=\tor_{M'}\text{ and }F_{(U,R)}\left(\tof_M\cap\W(U,R)\right)=\tof_{M'}$$
		where $(M',P')=E_{(U,R)}(M,P)$
	\end{corollary}
	
	\begin{proof}
		Let $E_{(U,R)}=F_{(U,R)}\circ\mathcal{E}_{(U,R)}$. By Proposition~\ref{l:redBM} and Lemma~\ref{l:wide}, this is a bijection from $\opname{\tau-rigid-pair}_{(U,R)}A$ to $\opname{\tau-rigid-pair} A_{(U,R)}$, commuting with direct sums and for any $(M,P)\in\opname{\tau-rigid-pair}_{(U,R)}A$, we have
		$\opname{co-rank}(M,P)=\opname{co-rank}E_{(U,R)}(M,P).$ 
		So $E_{(U,R)}$ restricts to a bijection from $\opname{s\tau-tilt-pair}_{(U,R)}A$ to $\opname{s\tau-tilt-pair} A_{(U,R)}$. 
		
		By Proposition~\ref{l:redBM} and Lemma~\ref{l:wide} again, for any $(M,P)\in \opname{s\tau-tilt-pair}_{(U,R)}A$, we have $F_{(U,R)}\left(\tor_M\cap\W(U,R)\right)=\tor_{M'}$, which implies that $E_{(U,R)}$ preserves the order. Since $\opname{s\tau-tilting-pair}_{(U,R)}A$ is an interval of $\opname{s\tau-tilting-pair}A$ by Lemma~\ref{l:AD} and the order can be induced by the mutation, the map $E_{(U,R)}$ commutes with the mutation.
		
		For any $X\in\W(U,R)$, since $(\tor_M,\tof_M)$ is a torsion pair in $\mod A$, there is an exact sequence
		$$0\to t_M(X)\to X\to f_M(X)\to 0$$
		with $t_M(X)\in\tor_M$ and $f_M(X)\in\tof_M$. Since $X\in U^\perp$ and $U^\perp$ is closed under submodules, we have $t_M(X)\in U^\perp$. On the other hand, $t_M(X)\in\tor_M\subseteq\tor_{T_{(U,R)}}$. So $t_M(X)\in\tor_{T_{(U,R)}}\cap U^\perp=\W(U,R)$. Since $\W(U,R)$ is a wide subcategory of $\mod A$ by Lemma~\ref{l:wide}, we have $f_M(X)\in\W(U,R)$ too. This shows that $(\tor_M\cap\W(U,R),\tof_M\cap\W(U,R))$ is a torsion pair in $\W(U,R)$. Then $(F_{(U,R)}\left(\tor_M\cap\W(U,R)\right),F_{(U,R)}\left(\tof_M\cap\W(U,R)\right))$ is a torsion pair in $\mod A_{(U,R)}$, which together with $F_{(U,R)}\left(\tor_M\cap\W(U,R)\right)=\tor_{M'}$, implies that $F_{(U,R)}\left(\tof_M\cap\W(U,R)\right)=\tof_{M'}$.
	\end{proof}
	
	The following result is useful (cf. \cite[Corollary 1.2]{BH}).
	
	\begin{corollary}\label{c:useful}
		Let $(U,R)$ be a basic $\tau$-rigid pair in $\mod A$ and $(N,Q)$ a basic $\tau$-rigid pair in $\mod A_{(U,R)}$. Then there exists a basic $\tau$-rigid pair $(N',Q')$ in $\mod A$ such that $A_{(N',Q')}\cong (A_{(U,R)})_{(N,Q)}$. 
	\end{corollary}
	
	\begin{proof}
		By Corollary~\ref{c:key}, there is $(N',Q')\in\opname{\tau-rigid-pair}_{(U,R)}A$ such that $E_{(U,R)}(N',Q')=(N,Q)$. Since $E_{(U,R)}$ commutes with direct sums and preserves the order, we have $E_{(U,R)}(T_{(N',Q')},Q')=(T_{(N,Q)},Q)$ and $E_{(U,R)}(S_{(N',Q')},L_{(N',Q')})=(S_{(N,Q)},L_{(N,Q)})$. Then we have $$F_{(U,R)}\left(\tor_{T_{(N',Q')}}\cap\W(U,R)\right)=\tor_{T_{(N,Q)}}\text{ and }F_{(U,R)}\left(\tof_{S_{(N',Q')}}\cap\W(U,R)\right)=\tof_{S_{(N,Q)}}.$$
		Since $\tor_{T_{(N',Q')}}\subseteq\tor_{T_{(U,R)}}$ and $\tof_{S_{(N',Q')}}\subseteq\tof_{S_{(U,R)}}$, we have 
		$$\tor_{T_{(N',Q')}}\cap\tof_{S_{(N',Q')}}\subseteq\tor_{T_{(U,R)}}\cap\tof_{S_{(U,R)}}=\W(U,R),$$
		which implies
		$$\left(\tor_{T_{(N',Q')}}\cap\W(U,R)\right)\cap\left(\tof_{S_{(N',Q')}}\cap\W(U,R)\right)=\tor_{T_{(N',Q')}}\cap\tof_{S_{(N',Q')}}=\W(N',Q').$$
		Hence we have $F_{(U,R)}(\W(N',Q'))=\W(N,Q)$. By Lemma~\ref{l:wide}, we have $A_{(N',Q')}\cong (A_{(U,R)})_{(N,Q)}$ as required.
	\end{proof}

\end{document}